\documentclass[10pt, en]{article}
\usepackage[margin=1.5cm]{geometry}                		% See geometry.pdf to learn the layout options. There are lots.
\DeclareMathSizes{10}{9}{7}{5}

\usepackage[table]{xcolor} 		

\usepackage{graphicx}	
\usepackage{authblk}
\usepackage{abstract}

\usepackage{amssymb, amsthm, bbold, amsmath} 
\usepackage{booktabs}
\usepackage{enumitem}
\usepackage{natbib}
\usepackage{algorithm}
\usepackage[noend]{algpseudocode}

\usepackage[font=large]{subcaption}
\usepackage{colortbl}

\newcommand{\bx}{\bold{x}}

\newcommand{\bL}{\bold{L}}
\newcommand{\bB}{\bold{B}}
\newcommand{\bR}{\bold{R}}
\newcommand{\bQ}{\bold{Q}}
\newcommand{\bleta}{\boldsymbol{\eta}}
\newcommand{\by}{\bold{y}}
\newcommand{\bH}{\bold{H}}
\newcommand{\bM}{\bold{M}}
\newcommand{\bI}{\bold{I}}
\newcommand{\bA}{\bold{A}}
\newcommand{\bb}{\bold{b}}
\newcommand{\bp}{\bold{p}}
\newcommand{\bP}{\bold{P}}
\newcommand{\bC}{\bold{C}}
\newcommand{\br}{\bold{r}}
\newcommand{\bS}{\bold{S}}
\newcommand{\bs}{\bold{s}}
\newcommand{\bV}{\bold{V}}
\newcommand{\bv}{\bold{v}}
\newcommand{\bU}{\bold{U}}
\newcommand{\bu}{\bold{u}}
\newcommand{\bZ}{\bold{Z}}
\newcommand{\bz}{\bold{z}}
\newcommand{\bK}{\bold{K}}
\newcommand{\bW}{\bold{W}}
\newcommand{\bw}{\bold{w}}
\newcommand{\bT}{\bold{T}}
\newcommand{\blf}{\bold{f}}
\newcommand{\bF}{\bold{F}}
\newcommand{\bG}{\bold{G}}
\newcommand{\bY}{\bold{Y}}
\newcommand{\bE}{\bold{E}}

\newcommand{\bsTheta}{\boldsymbol{\Theta}}
\newcommand{\bsSigma}{\boldsymbol{\Sigma}}

\title{Randomised preconditioning for the forcing formulation of weak constraint 4D-Var}

\author[1]{\normalsize Ieva Dau\v{z}ickait\.{e}*}

\author[1,2]{\normalsize Amos S. Lawless}

\author[1,3]{\normalsize Jennifer A. Scott}

\author[1,4]{\normalsize Peter Jan van Leeuwen}

\affil[1]{\footnotesize School Of Mathematical, Physical and Computational Sciences, University of Reading,UK}

\affil[2]{ \footnotesize National Centre for Earth Observation, Reading, UK}

\affil[3]{\footnotesize Scientific Computing Department, STFC Rutherford Appleton Laboratory, UK}

\affil[4]{\footnotesize Department of Atmospheric Science, Colorado State University, USA}

\date{}

\begin{document}

\maketitle

\begin{abstract}
There is growing awareness that errors in the model equations cannot be ignored in data assimilation methods such as four-dimensional variational assimilation (4D-Var). If allowed for, more information can be extracted from observations, longer time windows are possible, and the minimisation process is easier, at least in principle. Weak constraint 4D-Var estimates the model error and minimises a series of linear least-squares cost functions, which can be achieved using the conjugate gradient (CG) method; minimising each cost function is called an inner loop. CG needs preconditioning to improve its performance. In previous work, limited memory preconditioners (LMPs) have been constructed using approximations of the eigenvalues and eigenvectors of the Hessian in the previous inner loop. If the Hessian changes significantly in consecutive inner loops, the LMP may be of limited usefulness. To circumvent this, we propose using randomised methods for low rank eigenvalue decomposition and use these approximations to cheaply construct LMPs using information from the current inner loop. Three randomised methods are compared. Numerical experiments in idealized systems show that the resulting LMPs perform better than the existing LMPs. Using these methods may allow more efficient and robust implementations of incremental weak constraint 4D-Var.
\vspace{\baselineskip}

\end{abstract}

\let\thefootnote\relax\footnotetext{*Correspondence: Ieva Dau\v{z}ickait\.{e}, Department of Mathematics and Statistics, Pepper Lane, Whiteknights, Reading RG6 6AX, UK. Email: i.dauzickaite@pgr.reading.ac.uk}

\section{Introduction}
In numerical weather prediction, data assimilation provides the initial conditions for the weather model and hence influences the accuracy of the forecast \citep{Kalnay2002}. Data assimilation uses observations of a dynamical system to correct a previous estimate (background) of the system's state. The statistical knowledge of the errors in the observations and the background is incorporated in the process. A variational data assimilation method called weak constraint 4D-Var provides a way to also take into account the model error \citep{Tremolet06}, which can lead to a better analysis (e.g. \cite{Tremolet07}).

We explore the weak constraint 4D-Var cost function. In its incremental version, the state is updated by a minimiser of the linearised version of the cost function. The minimiser can be found by solving a large sparse linear system. %with the Hessian matrix 
The process of solving each system is called an inner loop. Because the second derivative of the cost function, the Hessian, is symmetric positive definite, the systems may be solved with the conjugate gradient (CG) method \citep{Hestenes52}, whose convergence rate depends on the eigenvalue distribution of the Hessian. Limited memory preconditioners (LMPs) have been shown to improve the convergence of CG when minimising the strong constraint 4D-Var cost function \citep{Fisher98, Tshimanga08}. Strong constraint 4D-Var differs from the weak constraint 4D-Var by making the assumption that the dynamical model has no error.

LMPs can be constructed using approximations to the eigenvalues and eigenvectors (eigenpairs) of the Hessian. The Lanczos and CG connection (Section 6.7 of \cite{SaadBook}) can be exploited to obtain approximations to the eigenpairs of the Hessian in one inner loop, and these approximations may then be employed to construct the preconditioner for the next inner loop \citep{Tshimanga08}. This approach does not describe how to precondition the first inner loop, and the number of CG iterations used on the $i$th inner loop limits the number of vectors available to construct the preconditioner on the $(i+1)$st inner loop. Furthermore, the success of preconditioning relies on the assumption that the Hessians do not change significantly from one inner loop to the next.

In this paper, we propose addressing these drawbacks by using the easy to implement subspace iteration methods (see Chapter 5 of \cite{SaadLEP}) to obtain approximations of the largest eigenvalues and corresponding eigenvectors of the Hessian in the current inner loop. The subspace iteration method first approximates the range of the Hessian by multiplying it with a start matrix (for approaches to choosing it see e.g. \cite{Duff93}) and the speed of convergence depends on the choice of this matrix (e.g. \cite{Gu2015}). A variant of subspace iteration, which uses a Gaussian random start matrix, is called the Randomised Eigenvalue Decomposition (REVD). REVD has been popularised by probabilistic analysis \citep{Halko11, Martinsson20}. %Unlike the usual subspace iteration method, REVD  to approximate the range of the Hessian. 
% Interest in this method that uses a Gaussian random matrix to approximate the range of the Hessian has been revived by probabilistic analysis \citep{Halko11, Martinsson20}, where it is called the Randomised Eigenvalue Decomposition (REVD). 
It has been shown that REVD, which is equivalent to one iteration of the subspace iteration method, can often generate a satisfactory approximation of the largest eigenpairs of a matrix that has rapidly decreasing eigenvalues. Because the Hessian is symmetric positive definite, a randomised Nystr\"{o}m method for computing a low rank eigenvalue decomposition can also be used. It is expected to give a higher quality approximation than REVD (e.g. \cite{Halko11}). We explore these two methods and another implementation of REVD, which is based on the \textit{ritzit} implementation of the subspace method \citep{Rutishauser1971}. The methods differ in the number of matrix-matrix products with the Hessian. Even though more computations are required to generate the preconditioner in the current inner loop compared to using information from the previous inner loop, the randomised methods are block methods and hence easily parallelisable.  %The comparison of the three methods show that they have different sensitivities, and one may be more suitable for data assimilation. We show that fewer vectors can be used in the preconditioner if they are obtained in the current inner loop with the randomised methods rather than in the previous inner loop using the Lanczos connection. 

In Section 2, we discuss the weak constraint 4D-Var method and, in Section 3, we consider LMPs and ways to obtain spectral approximations. The three randomised methods are examined in Section 4. Numerical experiments with linear advection and Lorenz 96 models are presented in Section 5, followed by a concluding discussion in Section 6.

\section{Weak constraint 4D-Var}
We are interested in estimating the state evolution of a dynamical system $\bx_0, \bx_1,\dots,\bx_N$, with $\bx_i \in \mathbb{R}^n$, at times $t_0, t_1,\dots,t_N$. Prior information about the state at $t_0$ is called the background and is denoted by $\bx^b \in \mathbb{R}^n$. It is assumed that $\bx^b$ has Gaussian errors with zero mean and covariance matrix $\bB \in \mathbb{R}^{n \times n}$. Observations of the system at time $t_i$ are denoted by $\by_i \in \mathbb{R}^{q_i}$ and their errors are assumed to be Gaussian with zero mean and covariance matrix $\bR_i \in \mathbb{R}^{q_i \times q_i}$ ($q_i \ll n$). An observation operator $\mathcal{H}_i$ maps the model variables into the observed quantities at the correct location, i.e. $\by_i = \mathcal{H}_i (\bx_i) + \boldsymbol{\zeta}_i$, where $\boldsymbol{\zeta}_i$ is the observation error. We assume that the observation errors are uncorrelated in time. 

The dynamics of the system are described using a nonlinear model $\mathcal{M}_i$ such that
\begin{equation}\label{eq:model}
\bx_{i+1} = \mathcal{M}_i (\bx_i) + \bleta_{i+1},
\end{equation}
where $\bleta_{i+1}$ is the model error at time $t_{i+1}$. The model errors are assumed to be uncorrelated in time and to be Gaussian with zero mean and covariance matrix $\bQ_i \in \mathbb{R}^{n \times n}$.

The forcing formulation of the nonlinear weak constraint 4D-Var cost function in which we solve for the initial state and the model error realizations, is\par
\begin{align}
J(\bx_0, \bleta_1, \dots, \bleta_N) & =  \frac{1}{2} (\bx_0 - \bx^b)^T \bB^{-1} (\bx_0 - \bx^b) + \frac{1}{2} \sum_{i=0}^{N} (\by_i - \mathcal{H}_i (\bx_i))^T \bR_i^{-1} (\by_i - \mathcal{H}_i (\bx_i)) \label{eq:4D-var_error} \\ \nonumber
& + \frac{1}{2} \sum_{i=1}^{N} \bleta_i^T \bQ_i^{-1} \bleta_i,
\end{align}
where $\bx_i$ satisfies the model equation \eqref{eq:model} \citep{Tremolet06}. The analysis (approximation of the state evolution over the time window) $\bx^a_0, \bx^a_1,\dots,\bx^a_N$ can be obtained from the minimiser of \eqref{eq:4D-var_error} using constraints \eqref{eq:model}.

\subsection{Incremental 4D-Var}
One way to compute the analysis is to approximate the minimum of \eqref{eq:4D-var_error} with an inexact Gauss-Newton algorithm \citep{Gratton2007}, where a sequence of quadratic cost functions is minimised. In this approach, we update $\bx_0$ and the model error  
\begin{equation}\label{eq:error_update}
\bold{p}^{(j+1)}  = \bold{p}^{(j)} +\delta \bold{p}^{(j)},  
\end{equation}
where $\bold{p}^{(j)} = (\bx_0^{(j)T},\bleta_1^{(j)T},\dots, \bleta_N^{(j)T})^T$ is the $j$th approximation and $\delta \bold{p}^{(j)} = (\delta \bx_0^{(j)T}, \delta \bleta_1^{(j)T},\dots, \delta \bleta_N^{(j)T})^T$. The $j$th approximation of the state $\bx^{(j)}= (\bx_0^{(j)\ T},\dots, \bx_N^{(j)\ T})^T$ is calculated with \eqref{eq:model} using $\bold{p}^{(j)}$. The update $\delta \bold{p}^{(j)}$ is obtained by minimising the following cost function %\delete{(we use bold notation for matrices and vectors that include the time dimension)}
\begin{equation}\label{eq:incr_wc4d-var_forcing}
J^{\delta} (\delta \bold{p}^{(j)}) = \frac{1}{2} || \delta \bold{p}^{(j)}- \bold{b}^{(j)} ||^2_{\bold{D}^{-1}} + \frac{1}{2} || \bold{H}^{(j)} (\bold{L}^{-1})^{(j)} \delta \bold{p}^{(j)} - \bold{d}^{(j)} ||^2_{\bold{R}^{-1}},
\end{equation}
where $||\bold{a}||^2_{\bold{A}^{-1}}=\bold{a}^T\bold{A}^{-1}\bold{a}$ and the covariance matrices are block diagonal, i.e. $\bold{D} = diag(\bB,\bQ_1,\dots,\bQ_N) \in \mathbb{R}^{n(N+1) \times n(N+1)}$ and $\bold{R} = diag(\bR_0, \dots, \bR_N) \in \mathbb{R}^{q \times q} $, where $q = \Sigma_{i=0}^{N} q_i$. We use the notation (following \cite{Gratton18}) $\bold{H}^{(j)} = diag(\bH_0^{(j)}, \dots, \bH_N^{(j)}) \in \mathbb{R}^{q \times n(N+1) }$, where $\bH_i^{(j)}$ is the linearised observation operator, and 

\begin{align}
(\bL^{-1})^{(j)} = & \left( \begin{array}{ccccc}
\bI         &             &            &                    & \\
\bM_{0,0}^{(j)}    &     \bI       &            &                    &  \\
\bM_{0,1}^{(j)}   &   \bM_{1,1}^{(j)}   &       \bI     &                    &  \\
\vdots              &   \vdots        & \ddots & \ddots        &  \\
\bM_{0,N-1}^{(j)}   &   \bM_{1,N-1}^{(j)}          &      \cdots        & \bM_{N-1,N-1}^{(j)}    & \bI  \\
\end{array} \right), \
\\
 \bold{b}^{(j)} = & \left( \begin{array}{c}
\bx^b - \bx_0^{(j)}\\
- \bleta_1^{(j)} \\
\vdots \\
-\bleta_N^{(j)}
\end{array} \right), \ 
\\
\bold{d}^{(j)} = & \left( \begin{array}{c}
\by_0  - \mathcal{H}_0 (\bx_0^{(j)})\\
\by_1  - \mathcal{H}_1 (\bx_1^{(j)}) \\
\vdots \\
\by_N  - \mathcal{H}_N (\bx_N^{(j)}) 
\end{array} \right),
\end{align}
where $\bM_{i,l}^{(j)} =\bM_l^{(j)}  \dots \bM_i^{(j)}$ and $\bM_i^{(j)} $ is the linearised model, i.e. $\bM_{i,l}^{(j)}$ denotes the linearised model integration from time $t_i$ to $t_{l+1}$, $(\bL^{-1})^{(j)} \in \mathbb{R}^{n(N+1) \times n(N+1)}$, $\bx^{(j)}, \delta \bx^{(j)}, \bold{b}^{(j)} \in \mathbb{R}^{n(N+1)}$ and $\bold{d}^{(j)}  \in \mathbb{R}^q$. The outer loop consists of updating \eqref{eq:error_update}, calculating $\bx^{(j)}, \bold{b}^{(j)}, \bold{d}^{(j)}$, and linearising $\mathcal{H}_i$ and $\mathcal{M}_i$ for the next inner loop.
% Minimising \eqref{eq:incr_wc4d-var_forcing} is called the inner loop.

The minimum of the quadratic cost function \eqref{eq:incr_wc4d-var_forcing} can be found by solving a linear system

\begin{align}
\boldsymbol{\mathcal{A}}^{(j)} \delta \bold{p}^{(j)} & = \bold{D}^{-1} \bold{b}^{(j)}  + (\bL^{-T})^{(j)} (\bold{H}^T)^{(j)}  \bold{R}^{-1} \bold{d^{(j)}}, \label{eq:forcing_form} \\ 
\boldsymbol{\mathcal{A}}^{(j)} & =(\bold{D}^{-1} + (\bL^{-T})^{(j)} (\bold{H}^T)^{(j)}  \bold{R}^{-1} (\bold{H})^{(j)} (\bL^{-1})^{(j)} ) \in \mathbb{R}^{n(N+1) \times n(N+1)},  \label{eq:forcing_form_matrix}
\end{align}
where $\boldsymbol{\mathcal{A}}^{(j)}$ is the Hessian of \eqref{eq:incr_wc4d-var_forcing}, which is symmetric positive definite. These large sparse systems are usually solved with the conjugate gradient (CG) method, whose convergence properties depend on the spectrum of $\boldsymbol{\mathcal{A}}^{(j)}$ (see Section~\ref{sec:PCG} for a discussion). In general, clustered eigenvalues result in fast convergence. We consider a technique to cluster eigenvalues of $\boldsymbol{\mathcal{A}}^{(j)}$ in the following section. From now on we omit the superscript $(j)$.

\subsection{Control Variable Transform}
A control variable transform, which is also called first level preconditioning, maps the variables $\delta \bold{p}$ to $\delta \tilde{\bold{p} }$, whose errors are uncorrelated (see, e.g. Section 3.2 of \cite{Lawless2013}). This can be denoted as the transformation $\bold{D}^{1/2} \delta \tilde{\bold{p} } =\delta  \bold{p}$, where $\bold{D} = \bold{D}^{1/2} \bold{D}^{1/2}$ and $\bold{D}^{1/2}$ is the symmetric square root. The update $\delta \tilde{\bold{p} }$ is then the solution of
\begin{equation}\label{eq:forcing_1stlvl_prec}
\boldsymbol{\mathcal{A}}^{pr} \delta \tilde{\bold{p} }= \bold{D}^{-1/2} \bold{b}  + \bold{D}^{1/2} \bL^{-T} \bold{H}^T  \bold{R}^{-1} \bold{d}, \quad \boldsymbol{\mathcal{A}}^{pr}=\bold{I} + \bold{D}^{1/2} \bL^{-T} \bold{H}^T  \bold{R}^{-1} \bold{H} \bL^{-1} \bold{D}^{1/2}.
\end{equation} 
$\boldsymbol{\mathcal{A}}^{pr}$ is the sum of the identity matrix and a rank $q$ positive semi-definite matrix. Hence, it has a cluster of $n(N+1) - q$ eigenvalues at one and $q$ eigenvalues that are greater than one. Thus, the spectral condition number $\kappa = \lambda_{max} / \lambda_{min}$ (here $\lambda_{max}$ and $\lambda_{min} $ are the largest and smallest, respectively, eigenvalues of $\boldsymbol{\mathcal{A}}^{pr}$) is equal to $\lambda_{max}$. We discuss employing second level preconditioning to reduce the condition number while also preserving the cluster of the eigenvalues at one. In the subsequent sections, we use notation that is common in numerical linear algebra. Namely, we use $\bA$ for the Hessian with first level preconditioning, $\bx$ for the unknown and $\bb$ for the right hand side of the system of linear equations. Thus, we denote \eqref{eq:forcing_1stlvl_prec} by 
\begin{equation}\label{eq:Ax_eq_b}
\bA \bx=\bb,
\end{equation}
where the right-hand side $\bb$ is known and $\bx$ is the required solution. We assume throughout that $\bA$ is symmetric positive definite.

\section{Preconditioning weak constraint 4D-Var}

\subsection{Preconditioned conjugate gradients}\label{sec:PCG}

The CG method (see, e.g. \cite{SaadBook}) is a popular Krylov subspace method for solving systems of the form \eqref{eq:Ax_eq_b}. A well known bound for the error at the $i$th CG iteration $\boldsymbol{\epsilon}_i = \bx - \bx_i$  is
\begin{equation}
\frac{|| \boldsymbol{\epsilon}_i ||_\bA }{|| \boldsymbol{\epsilon}_0 ||_\bA} \leq 2 \left( \frac{\sqrt{\kappa} - 1}{\sqrt{\kappa} + 1} \right)^i,
\end{equation}
where $\kappa$ is the spectral condition number and $|| \boldsymbol{\epsilon}_i ||^2_\bA = \boldsymbol{\epsilon}_i^T \bA \boldsymbol{\epsilon}_i$ (see, e.g., Section 5.1. of \cite{Nocedal06}). Note that this bound describes the worst-case convergence and only takes into account the largest and smallest eigenvalues. The convergence of CG also depends on the distribution of the eigenvalues of $\bA$ (as well as the right hand side $\bb$ and the initial guess $\bx_0$); eigenvalues clustered away from zero suggest rapid convergence (Lecture 38 of \cite{TrefethenBau}). Otherwise, CG can display slow convergence and preconditioning is used to try and tackle this problem (Chapter 9 of \cite{SaadBook}). Preconditioning aims to map the system \eqref{eq:Ax_eq_b} to another system that has the same solution, but different properties that imply faster convergence. Ideally, the preconditioner $\bP$ should be both cheap to construct and to apply, and the preconditioned system should be easy to solve.

If $\bP$ is a symmetric positive definite matrix that approximates $\bA^{-1}$ and is available in factored form $\bP = \bC \bC^T$, the following system is solved
\begin{equation}\label{eq:split_prec_system}
\bC^T \bA \bC \hat{\bx} = \bC^T \bb,
\end{equation}
where $\hat{\bx} = \bC^{-1} \bx$. Split preconditioned CG (PCG) for solving \eqref{eq:split_prec_system} is described in Algorithm~\ref{alg:PCG} (see, for example, Algorithm 9.2 of \cite{SaadBook}). Note that every CG iteration involves one matrix-vector product with $\bA$ (the product $\bA \bp_{j-1}$ is stored in step~\ref{step:alpha} and reused in step~\ref{step:res}) and this is expensive in weak constraint 4D-Var, because it involves running the linearised model throughout the assimilation window through the factor $\bL^{-1}$.

\begin{algorithm}
\caption{Split preconditioned CG for solving $\bA \bx=\bb$}\label{alg:PCG}
\hspace*{\algorithmicindent} \textbf{Input:} $\bA \in \mathbb{R}^{n_A \times n_A}$, $\bb \in \mathbb{R}^{n_A}$, preconditioner $\bP = \bC \bC^T \in \mathbb{R}^{n_A \times n_A}$, initial solution  $\bx_0 \in \mathbb{R}^{n_A}$ \\
 \hspace*{\algorithmicindent} \textbf{Output:} solution $\bx_j \in \mathbb{R}^{n_A}$
\begin{algorithmic}[1]
\State Compute  $\br_0 =\bC^T (\bb - \bA \bx_0)$, and $\bp_0 =\bC \br_0$
\For {$j=1,2, \dots,$ until convergence}
\State $\alpha_j = (\br_{j-1}^T \br_{j-1}) / (\bp_{j-1}^T \bA \bp_{j-1})$ \label{step:alpha}
\State $\bx_j = \bx_{j-1} + \alpha_j \bp_{j-1} $
\State $\br_j = \br_{j-1} - \alpha_j \bC^T \bA \bp_{j-1}$ \label{step:res}
\State $\beta_j = (\br_j^T \br_j) / (\br_{j-1}^T \br_{j-1})$
\State $\bp_j = \bC \br_j + \beta_j \bp_{j-1}$
\EndFor
\end{algorithmic}
\end{algorithm}

\subsection{Limited memory preconditioners}\label{sec:LMPs}
In weak constraint 4D-Var the preconditioner $\bP$ approximates the inverse Hessian. Hence, $\bP$ can be obtained using Quasi-Newton methods for unconstrained optimization that construct an approximation of the Hessian matrix, which is updated regularly (see, for example, Chapter 6 of \cite{Nocedal06}). A popular method to approximate the Hessian is BFGS (named after Broyden, Fletcher, Goldfarb, and Shanno, who proposed it), but it is too expensive in terms of storage and updating the approximation. Instead, the so-called block BFGS method (derived by \cite{Schnabel83}) uses only a limited number of vectors to build the Hessian, and when new vectors are added older ones are dropped. This is an example of a limited memory preconditioner (LMP), and the one considered by 
%A popular quasi-Newton method BFGS approximates the inverse of the Hessian. A generalisation of BFGS, called block BFGS (derived by \cite{Schnabel83}), is used to obtain the limited-memory preconditioners (LMPs) considered by 
Tshimanga et al. (see \cite{Tshimanga08, Gratton11} and \cite{Tshimanga_thesis}) in the context of strong constraint 4D-Var. An LMP for an $n_A \times n_A$ symmetric positive definite matrix $\bA$ is defined as follows
\begin{equation}\label{eq:LMP_general}
\bP_k = ( \bI_{n_A} - \bS (\bS^T \bA \bS)^{-1} \bS^T \bA) ( \bI_{n_A} - \bA \bS(\bS^T \bA \bS)^{-1} \bS^T) + \bS (\bS^T \bA \bS)^{-1} \bS^T,
\end{equation}
where $\bS$ is an $n_A \times k$ ($k \leq n_A$) matrix with linearly independent columns $\bs_1, \dots, \bs_k$, and $\bI_{n_A}$ is the $n_A \times n_A$ identity matrix \citep{Gratton11}. $\bP_k$ is symmetric positive definite and if $k = n_A$ then $(\bS^T \bA \bS)^{-1} = \bS^{-1} \bA^{-1} \bS^{-T}$ and $\bP_k = \bA^{-1}$. In data assimilation, we have $k \ll n_A$, hence the name LMPs. $\bP_k$ is called a balancing preconditioner in \cite{Tang09}.

A potential problem for practical applications of \eqref{eq:LMP_general} is the need for expensive matrix-matrix products with $\bA$. Simpler formulations of \eqref{eq:LMP_general} are obtained by imposing more conditions on the vectors $\bs_1, \dots, \bs_k$. Two formulations that \cite{Tshimanga08} calls spectral-LMP and Ritz-LMP have been used, for example, in ocean data assimilation in the Regional Ocean Modeling System (ROMS) \citep{Moore11} and the variational data assimilation software with the Nucleus for European Modelling of the Ocean (NEMO) ocean model (NEMOVAR) \citep{Nemovar_ecmwf}, and coupled climate reanalysis in Coupled ECMWF ReAnalysis (CERA) \citep{Laloyaux18}. 

To obtain the spectral-LMP, let $\bv_1, \dots, \bv_k$ be orthonormal eigenvectors of $\bA$ with corresponding eigenvalues $\lambda_1, \dots, \lambda_k$. Set $\bV = (\bv_1, \dots, \bv_k)$ and $\boldsymbol{\Lambda} = diag(\lambda_1, \dots, \lambda_k)$ so that $\bA \bV = \bV \boldsymbol{\Lambda}$ and $\bV^T \bV = \bI_k$. If $\bs_i = \bv_i$, $i=1,\dots,k$, then the LMP in \eqref{eq:LMP_general} is the spectral-LMP $\bP_k^{sp}$ (it is called a deflation preconditioner in \cite{Giraud06}), which can be simplified as
\begin{equation}\label{eq:spectral_LMP}
\bP_k^{sp} = \bI_{n_A} - \sum_{i=1}^{k} (1 - \lambda_i^{-1}) \bv_i \bv_i^T.
\end{equation}
Then $\bP_k^{sp} = \bC_k^{sp}(\bC_k^{sp})^T$ with (presented in Section 2.3.1 of \cite{Tshimanga_thesis})
\begin{equation}\label{eq:spectral_LMP_factor}
\bC_k^{sp} = \prod_{i=1}^{k} \left(\bI_{n_A} - \left( 1 - (\sqrt{\lambda_i})^{-1}\right) \bv_i \bv_i^T \right).
\end{equation}

In many applications, including data assimilation, exact eigenpairs are not known, and their approximations, called Ritz values and vectors, are used (we discuss these in the following section). %By setting $s_1, \dots, s_{k}$ to be orthogonal Ritz vectors $u_1, \dots, u_{k}$ and using the relation $U^T A U = \Theta$, where $U$ has columns $u_1, \dots, u_{k}$, $\Theta = diag(\theta_1, \dots, \theta_{k})$ and $\theta_i$ is a Ritz value, the Ritz-LMP $P_k^{Rt}$ is
If $\bu_1, \dots, \bu_k$ are orthogonal Ritz vectors, then the following relation holds $\bU^T \bA \bU = \boldsymbol{\Theta}$, where $\bU = (\bu_1, \dots, \bu_k)$, $\boldsymbol{\Theta} = diag(\theta_1, \dots, \theta_k)$ and $\theta_i$ is a Ritz value. Setting $\bs_i = \bv_i$, $i=1,\dots,k$, the Ritz-LMP $\bP_k^{Rt}$ is
\begin{equation}\label{eq:Ritz_LMP}
\bP_k^{Rt} =  ( \bI_{n_A} - \bU \boldsymbol{\Theta}^{-1} \bU^T \bA) ( \bI_{n_A} - \bA \bU \boldsymbol{\Theta}^{-1} \bU^T) + \bU \boldsymbol{\Theta}^{-1} \bU^T.
\end{equation}
Each application of $\bP_k^{Rt}$ requires a matrix-matrix product with $\bA$. If the Ritz vectors are obtained by the Lanczos process (described in Section~\ref{sec:Lanczos_process} below), then \eqref{eq:Ritz_LMP} can be further simplified, so that no matrix-matrix products with $\bA$ are needed (see Section 4.2.2. of \cite{Gratton11} for details).

%Then the factor $C_k^{Rt}$ of the Ritz-LMP is:
%\begin{equation}\label{eq:Ritz_LMP_factor}
%C_k^{Rt} = \prod_{i=1}^{k} \left( I_{n_A} - \left( 1 - (\sqrt{\theta_i})^{-1} \right) u_i u_i^T - \frac{e^T_k w_i}{\sqrt{\theta_i}} \tau_{k-1} u_i f_{k+1} ^T \right),
%\end{equation}
%where $w_i$ is a primitive Ritz vector, $e_k$ is a zero vector with 1 as its $k$th entry, $f_{k+1}$ is a Lanczos vector and $ \tau_{k-1}$ is an off-diagonal entry of matrix $T_k$ used in the Lanczos process. 
% split formulations are given in sections 4.2.1. and 4.2.2. of Gratton et al.

An important property is that if an LMP is constructed using $k$ vectors then at least $k$ eigenvalues of the preconditioned matrix $\bC^T \bA \bC$ will be equal to 1, and the remaining eigenvalues will lie between the smallest and largest eigenvalues of $\bA$ (see Theorem 3.4 of \cite{Gratton11}). Moreover, if $\bA$ has a cluster of eigenvalues at 1, then LMPs preserve this cluster. This is crucial when preconditioning \eqref{eq:forcing_1stlvl_prec}: because the LMPs preserve the $n(N+1) - q$ smallest eigenvalues of $\boldsymbol{\mathcal{A}}^{pr}$ that are equal to $1$, the CG convergence can be improved by decreasing the largest eigenvalues. Hence, it is preferable to use the largest eigenpairs or their approximations. 
%The latter is crucial when preconditioning \eqref{eq:forcing_1stlvl_prec}, because as already noted, the $n(N+1) - q$ smallest eigenvalues of $\mathcal{A}^{pr}$ are equal to $1$. In this case, the CG convergence can be adversely affected by the largest eigenvalues of $\mathcal{A}^{pr}$. Hence, it is preferable to use the eigenvectors and Ritz vectors associated with the largest eigenvalues and Ritz values. % Also Th. 3.7 talks about the multiplicity of the eigenvalue 1 when there already are some eigenvalues at 1. Talks about the cases with eigenvectors in the preconditioner.

In practice, both the spectral-LMP and Ritz-LMP use Ritz vectors and values to construct the LMPs. %Notice that the Ritz-LMP is more general, and $C_k^{Rt}$ is equal to $C_k^{sp}$ with a correction term. 
It has been shown that the Ritz-LMP can perform better than the spectral-LMP in a strong constraint 4D-Var setting by correcting for the inaccuracies in the estimates of eigenpairs \citep{Tshimanga08}. However, \cite{Gratton11} (Theorem 4.5) have shown that if the preconditioners are constructed with Ritz vectors and values that have converged, then the spectral-LMP acts like the Ritz-LMP. 

\subsection{Ritz information}
Calculating or approximating all the eigenpairs of a large sparse matrix is impractical. Hence, only a subset is approximated to construct the LMPs. This is often done by extracting approximations from a subspace, and the Rayleigh-Ritz (RR) procedure is a popular method for doing this.

Assume that $\mathcal{Z} \subset \mathbb{R}^{n_A}$ is an invariant subspace of $\bA$, i.e. $\bA \bz \in \mathcal{Z}$ for every $\bz \in \mathcal{Z}$, and the columns of $\bZ \in \mathbb{R}^{n_A \times m}$, $m < n_A$, form an orthonormal basis for $\mathcal{Z}$. If $(\lambda, \hat{\by})$ is an eigenpair of $\bK = \bZ^T \bA \bZ \in \mathbb{R}^{m \times m}$, then $(\lambda, \bv)$, where $\bv = \bZ \hat{\by}$, is an eigenpair of $\bA$ (see, e.g. Theorem 1.2 in Chapter 4 of \cite{StewartMTX_ALGO}). Hence, eigenvalues of $\bA$ that lie in the subspace $\mathcal{Z}$ can be extracted by solving a small eigenvalue problem. 

However, generally the computed subspace $\tilde{\mathcal{Z}}$ with orthonormal basis as columns of $\tilde{\bZ}$ is not invariant. Hence, only approximations $\tilde{\bv}$ to the eigenvectors $\bv$ belong to $\tilde{\mathcal{Z}}$. The RR procedure computes approximations $\bu$ to $\tilde{\bv}$. We give the RR procedure in Algorithm~\ref{alg:rr}, where the eigenvalue decomposition is abbreviated as EVD. Approximations to eigenvalues $\lambda$ are called Ritz values $\theta$, and $\bu$ are the Ritz vectors. Eigenvectors of $\tilde{\bK}=\tilde{\bZ}^T \bA \tilde{\bZ}$, which is the projection of $\bA$ onto $\tilde{\mathcal{Z}}$, are denoted by $\bw$ and are called primitive Ritz vectors.

\begin{algorithm}
\caption{Rayleigh-Ritz procedure for computing approximations of eigenpairs of symmetric $\bA$}\label{alg:rr}
\hspace*{\algorithmicindent} \textbf{Input:} symmetric matrix $\bA \in \mathbb{R}^{n_A \times n_A}$, orthogonal matrix $\tilde{\bZ} \in \mathbb{R}^{n_A \times m}$, $m<n_A$ \\
 \hspace*{\algorithmicindent} \textbf{Output:} orthogonal $\bU \in \mathbb{R}^{n_A \times m}$ with approximations to eigenvectors of $\bA$ as its columns, and diagonal $\boldsymbol{\Theta} \in \mathbb{R}^{m \times m}$ with approximations to eigenvalues of $\bA$ on the diagonal
\begin{algorithmic}[1]
\State  Form $\tilde{\bK}=\tilde{\bZ}^T \bA \tilde{\bZ} \in \mathbb{R}^{m \times m}$
\State Form EVD of $\tilde{\bK}:\ \tilde{\bK} = \bW \boldsymbol{\Theta} \bW^T$, where $\bW, \ \boldsymbol{\Theta} \in \mathbb{R}^{ m \times m}$
\State Form Ritz vectors $\bU = \tilde{\bZ} \bW \in \mathbb{R}^{n_A \times m}$
\end{algorithmic}
\end{algorithm}

\subsection{Spectral information from CG}\label{sec:Lanczos_process}
\cite{Tshimanga08} use Ritz pairs of the Hessian in one inner loop to construct LMPs for the following inner loop, i.e. information on $\bA^{(0)}$ is used to precondition $\bA^{(1)}$, and so on. Success relies on the Hessians not changing significantly from one inner loop to the next. Ritz information can be obtained from the Lanczos process that is connected to CG, hence information for the preconditioner can be gathered at a negligible cost.

The Lanczos process is used to obtain estimates of a few extremal eigenvalues and corresponding eigenvectors of a symmetric matrix $\bA$ (Section 10.1 of \cite{Golub13}). It produces a sequence of tridiagonal matrices $\bT_j \in \mathbb{R}^{j \times j}$, whose largest and smallest eigenvalues converge to the largest and smallest eigenvalues of $\bA$. Given a starting vector $\blf_0$, it also computes an orthonormal basis $\blf_0,\dots, \blf_{j-1}$ for the Krylov subspace $\mathcal{K}_j = span\{\blf_0, \bA \blf_0, \dots, \bA^{j-1} \blf_0\}$. %The residual $r_i = (A - \theta_i I )u_i$, where $(\theta_i, u_i)$ are the Ritz pairs of $\bA$, are orthogonal to $\mathcal{K}_j$.
Ritz values $\theta_i$ are obtained as eigenvalues of a tridiagonal matrix, which has the following structure:
\begin{equation}
\bT_j = \left( \begin{array}{cccc}
\gamma_1 & \tau_1 & 			&\\
\tau_1  & \gamma_2 & \tau_2	& \\
			& 	\ddots		& \ddots & \ddots \\
			&					&	\tau_{j-1}  & \gamma_j \\				
\end{array} \right).
\end{equation} 
The Ritz vectors of $\bA$ are $\bu_i=\bF_j \bw_i$, where $\bF_j=(\blf_0, \dots, \blf_{j-1})$ and an eigenvector $\bw_i$ of $\bT_j$ is a primitive Ritz vector. Eigenpairs of $\bT_j$ can be obtained using a symmetric tridiagonal QR algorithm or Jacobi procedures (e.g. Section 8.5 of \cite{Golub13}).

Saad (Section 6.7.3 of \cite{SaadBook}) discusses obtaining entries of $\bT_j$ when solving $\bA \bx=\bb$ with CG. At the $j$-th iteration of CG, new entries of $\bT_j$ are calculated as follows

\begin{gather}
\gamma_j = \begin{cases} 
      \frac{1}{\alpha_j} & \text{for } j =1 \\
      \frac{1}{\alpha_j} + \frac{\beta_{j-1}}{\alpha_{j-1}} & \text{for } j>1 \\
   \end{cases}\\
\tau_j = \frac{\sqrt{\beta_j}}{\alpha_j},
\end{gather}
and the vector $\blf_j = \br_j/ ||\br_j||$, where $||\br_j||^2 = \br_j^T \br_j$ and $\alpha_j, \beta_j$ and $\br_j$ are as in Algorithm~\ref{alg:PCG}. Hence, obtaining eigenvalue information requires normalizing the residual vectors and finding eigenpairs of the tridiagonal matrix $\bT_j$. In data assimilation, the dimension of $\bT_j$ is small, because the cost of matrix-vector products restricts the number of CG iterations in the previous inner loop. Hence its eigenpairs can be calculated cheaply. However, caution has to be taken to avoid `ghost' values, i.e. repeated Ritz values, due to the loss of orthogonality in CG (Section 10.3.5 of \cite{Golub13}). This can be addressed using a complete reorthogonalization in every CG iteration, which is done in the CONGRAD routine used at the European Centre for Medium Range Weather Forecasts \citep{IFSdocDA2020}. This makes every CG iteration more expensive, but CG may converge in fewer iterations \citep{Fisher98}.

\section{Randomised eigenvalue decomposition}\label{sec:randomised_evd}
%Although information for constructing the preconditioner can be approximated cheaply from the Lanczos process, in \cite{Tshimanga08} the LMP is constructed with eigenpairs of the Hessian from the previous inner loop. Each iteration of the Lanczos process requires a matrix-vector product with the Hessian, thus the cost is similar to the cost of CG and is too prohibitive for it to be used to obtain information on the current Hessian. Hence, we explore a different approach.
If the Hessian in one inner loop differs significantly from the Hessian in the previous inner loop, then it may not be useful to precondition the former with an LMP that is constructed with information from the latter. Employing the Lanczos process to obtain eigenpair estimates and use them to construct an LMP in the same inner loop is too computationally expensive, because each iteration of the Lanczos process requires a matrix-vector product with the Hessian, thus the cost is similar to the cost of CG. Hence, we explore a different approach.

Subspace iteration is a simple procedure to obtain approximations to the largest eigenpairs (see, e.g., Chapter 5 of \cite{SaadLEP}). It is easily understandable and can be implemented in a straightforward manner, although its convergence can be very slow if the largest eigenvalues are not well separated from the rest of the spectrum. The accuracy of subspace iteration may be enhanced by using a RR projection.

Such an approach is used in the Randomised Eigenvalue Decomposition (REVD) (see, e.g., \cite{Halko11}). This takes a Gaussian random matrix, i.e. a matrix with independent standard normal random variables with zero mean and variance equal to one as its entries, and applies one iteration of the subspace iteration method with RR projection, hence obtaining a rank $m$ approximation $\bA\approx \bZ_1 (\bZ_1^T \bA \bZ_1) \bZ_1^T$, where $\bZ_1 \in \mathbb{R}^{n_A \times m}$ is orthogonal. We present REVD in Algorithm~\ref{alg:revd}. An important feature of REVD is the observation that the accuracy of the approximation is enhanced with oversampling (which is also called using guard vectors in \cite{Duff93}), i.e. working on a larger space than the required number of Ritz vectors. \cite{Halko11} claim that setting the oversampling parameter to 5 or 10 is often sufficient.

\begin{algorithm}
\caption{Randomised eigenvalue decomposition, REVD}\label{alg:revd}
\hspace*{\algorithmicindent} \textbf{Input:} symmetric matrix $\bA \in \mathbb{R}^{n_A \times n_A}$,  target rank $k$, an oversampling parameter $l$ \\
 \hspace*{\algorithmicindent} \textbf{Output:} orthogonal $\bU_1 \in \mathbb{R}^{n_A \times k}$  with approximations to eigenvectors of $\bA$ as its columns, and diagonal $\bsTheta_1 \in \mathbb{R}^{k \times k}$ with approximations to the largest eigenvalues of $\bA$ on the diagonal
\begin{algorithmic}[1]
\State  Form a Gaussian random matrix $\bG \in \mathbb{R}^{n_A \times (k+l)}$
\State  Form a sample matrix $\bY=\bA \bG \in \mathbb{R}^{n_A \times (k+l)}$
\State  Orthonormalize the columns of $\bY$ to obtain orthonormal $\bZ_1 \in \mathbb{R}^{n_A \times (k+l)}$
\State  Form $\bK_1 = \bZ_1^T \bA \bZ_1 \in \mathbb{R}^{ (k+l) \times (k+l) }$ \label{alg_step:K=ZTAZ}
\State Form EVD of $\bK:\ \bK = \bW_1 \bsTheta_1 \bW_1^T$, where $\bW_1, \  \bsTheta_1 \in \mathbb{R}^{ (k+l) \times (k+l)}$, elements of $\bsTheta_1$ are sorted in decreasing order
\State Remove last $l$ columns and rows of $\bsTheta_1$, so that $\bsTheta_1 \in \mathbb{R}^{ k \times k}$
\State Remove last $l$ columns of $ \bW_1$, so that $\bW_1 \in \mathbb{R}^{ (k+l) \times k}$
\State Form $\bU_1 =  \bZ_1 \bW_1 \in \mathbb{R}^{n_A \times k}$.
\end{algorithmic}
\end{algorithm}

Randomised algorithms are designed to minimise the communication instead of the flop count. The expensive parts of Algorithm~\ref{alg:revd} are the two matrix-matrix products $\bA \bG$ and $\bA \bZ_1$ in steps 2 and 4, i.e. matrix $\bA$ has to be multiplied with $2(k+l)$ vectors, which in serial computations would be essentially the cost of $2(k+l)$ iterations of unpreconditioned CG. However, note that these matrix-matrix products can be parallelised. 

In weak constraint 4D-Var, $\bA$ is the Hessian, hence it is symmetric positive definite and its eigenpairs can also be approximated using a randomised Nystr\"{o}m method (Algorithm 5.5. of \cite{Halko11}), which is expected to give much more accurate results than REVD \citep{Halko11}. We present the Nystr\"{o}m method in Algorithm~\ref{alg:nystrom}, where singular value decomposition is abbreviated as SVD. It considers a more elaborate rank $m$ approximation than in REVD: $\bA\approx (\bA \bZ_1) (\bZ_1^T \bA \bZ_1)^{-1} (\bA \bZ_1)^T = \bF \bF^T$, where $\bZ_1 \in \mathbb{R}^{n_A \times m}$ is orthogonal (obtained in the same way as in REVD, e.g. using a tall skinny QR (TSQR) decomposition\citep{Demmel2012}) and $\bF =(\bA \bZ_1) (\bZ_1^T \bA \bZ_1)^{-1/2} \in \mathbb{R}^{n_A \times m}$ is an approximate Cholesky factor of $\bA$, which is found in step 6. The eigenvalues of $\bF \bF^T$ are the squares of the singular values of $\bF$ (see section 2.4.2 of \cite{Golub13}). In numerical computations we store matrices $\bE^{(1)} = \bA \bZ_1$ and $\bE^{(2)} =  \bZ_1^T \bE^{(1)} = \bZ_1^T \bA \bZ_1$ (step 4), perform the Cholesky factorization of $\bE^{(2)} =  \bC^T \bC$ (step 5) and obtain $\bF$ by solving the triangular system $\bF \bC= \bE^{(1)}$. 

\begin{algorithm}
\caption{Randomised eigenvalue decomposition for symmetric positive semidefinite $\bA$, Nystr\"{o}m}\label{alg:nystrom}
\hspace*{\algorithmicindent} \textbf{Input:} symmetric positive semidefinite matrix $\bA \in \mathbb{R}^{n_A \times n_A}$,  target rank $k$, an oversampling parameter $l$ \\
 \hspace*{\algorithmicindent} \textbf{Output:} orthogonal $ \bU_2 \in \mathbb{R}^{n_A \times k}$  with approximations to eigenvectors of $\bA$ as its columns, and diagonal $ \bsTheta_2 \in \mathbb{R}^{k \times k}$ with approximations to eigenvalues of $\bA$ on the diagonal
\begin{algorithmic}[1]
\State  Form a Gaussian random matrix $\bG \in \mathbb{R}^{n_A \times (k+l)}$
\State  Form a sample matrix $\bY=\bA \bG \in \mathbb{R}^{n_A \times (k+l)}$
\State  Orthonormalize the columns of $\bY$ to obtain orthonormal $ \bZ_1 \in \mathbb{R}^{n_A \times (k+l)}$
\State  Form matrices $\bE^{(1)} = \bA \bZ_1 \in \mathbb{R}^{n_A \times (k+l)}$ and $\bE^{(2)} = \bZ_1^T \bE^{(1)} \in \mathbb{R}^{ (k+l) \times (k+l) }$ 
\State Form a Cholesky factorization $\bE^{(2)} = \bC^T \bC$
\State Solve $\bF \bC= \bE^{(1)}$ for $\bF \in \mathbb{R}^{n_A \times (k+l)}$
\State Form SVD of $\bF:\ \bF = \bU_2 \bsSigma \bV^T $, where $\bU_2, \ \bV \in \mathbb{R}^{n_A\times (k+l)}$, $\bsSigma \in \mathbb{R}^{ (k+l) \times (k+l)}$, elements of $\bsSigma$ are sorted in decreasing order
\State Remove last $l$ columns of $\bU_2$, so that $\bU_2 \in \mathbb{R}^{n_A \times k}$
\State Remove last $l$ columns and rows of $\bsSigma$, so that $\bsSigma \in \mathbb{R}^{ k \times k}$, and set $\bsTheta_2 = \bsSigma^2$
\end{algorithmic}
\end{algorithm}

The matrix-matrix product with $\bA$ at step \ref{alg_step:K=ZTAZ} of Algorithms~\ref{alg:revd} and \ref{alg:nystrom} is removed in Rutishauser's implementation of subspace iteration with RR projection called \textit{ritzit} \citep{Rutishauser1971}. It can be derived in the following manner (see Chapter 14 of \cite{Parlett_TSEP}). Assume that $\bG_3 \in \mathbb{R}^{n_A \times m}$ is an orthogonal matrix and the sample matrix is $\bY_3 = \bA \bG_3 = \bZ_3 \bR_3$, where $\bZ_3 \in \mathbb{R}^{n_A \times m}$ is orthogonal and $\bR_3 \in \mathbb{R}^{m \times m}$ is upper triangular. Then a projection of $\bA^2$ onto the column space of $\bG_3$ is $\hat{\bK} = \bY_3^T \bY_3 = \bR_3^T \bZ_3^T \bZ_3 \bR_3 = \bR_3^T \bR_3$. Then $\bK_3 = \bR_3 \bR_3^T = \bR_3 \bR_3^T \bR_3 \bR_3^{-1} = \bR_3 \hat{\bK} \bR_3^{-1}$, which is similar to $\hat{\bK}$ and hence has the same eigenvalues. This leads to another implementation of REVD presented in Algorithm~\ref{alg:ritzit}. This is a single pass algorithm, meaning that $\bA$ has to be accessed just once, and to the best of our knowledge this method has not been considered in the context of randomised eigenvalue approximations.

\begin{algorithm}
\caption{Randomised eigenvalue decomposition based on \textit{ritzit}, REVD\_\textit{ritzit} }\label{alg:ritzit}
\hspace*{\algorithmicindent} \textbf{Input:} symmetric matrix $\bA \in \mathbb{R}^{n_A \times n_A}$,  target rank $k$, an oversampling parameter $l$ \\
 \hspace*{\algorithmicindent} \textbf{Output:} orthogonal $\bU_3 \in \mathbb{R}^{n_A \times k}$ with approximations to eigenvectors of $\bA$ as its columns, and diagonal $\bsTheta_3 \in \mathbb{R}^{k \times k}$ with approximations to eigenvalues of $\bA$ on the diagonal
\begin{algorithmic}[1]
\State  Form a Gaussian random matrix $\bG \in \mathbb{R}^{n_A \times (k+l)}$
\State Orthonormalize the columns of $\bG$ to obtain orthonormal $\bG_3$
\State Form a sample matrix $\bY_3=\bA \bG_3 \in \mathbb{R}^{n_A \times (k+l)}$
\State Compute QR decomposition $\bY_3=\bZ_3 \bR_3$ to obtain orthogonal $\bZ_3 \in \mathbb{R}^{n_A \times (k+l)}$ and upper triangular $\bR_3 \in \mathbb{R}^{(k+l) \times (k+l)}$
\State Form $\bK_3 =  \bR_3  \bR_3^T \in \mathbb{R}^{ (k+l) \times (k+l) }$
\State Form EVD of $\bK_3:\ \bK_3 = \bW_3 \bsTheta_3^2  \bW_3^T$, where $\bW_3, \ \bsTheta_3^2 \in \mathbb{R}^{ (k+l) \times (k+l)}$, elements of $\bsTheta_3$ are sorted in decreasing order
\State Remove last $l$ columns and rows of $\bsTheta_3^2$, so that $\bsTheta_3^2 \in \mathbb{R}^{ k \times k}$
\State Remove last $l$ columns of $ \bW_3$, so that $ \bW_3 \in \mathbb{R}^{ (k+l) \times k}$
\State Form $\bU_3 = \bZ_3 \bW_3 \in \mathbb{R}^{n_A \times k}$.
\end{algorithmic}
\end{algorithm}

Note that the Ritz vectors given by Algorithms~\ref{alg:revd}, \ref{alg:nystrom} and \ref{alg:ritzit} are different. Although Algorithm~\ref{alg:ritzit} accesses the matrix $\bA$ only once, it requires an additional orthogonalisation of a matrix of size $n_A \times (k+l)$.

In Table~\ref{table:all_formulations_summary}, we summarise some properties of the Lanczos, REVD, Nystr\"{o}m and REVD\_\textit{ritzit}  methods when they are used to compute Ritz values and vectors to generate a preconditioner for $\bA$ in incremental data assimilation. Note that the cost of applying the spectral-LMP depends on the number of vectors $k$ used in its construction and is independent of which method is used to obtain them. The additional cost of using randomised algorithms arises only once per inner loop when the preconditioner is generated. We recall that in these algorithms the required EVD or SVD of the small matrix can be obtained cheaply and the most expensive parts of the algorithms are the matrix-matrix products of $\bA$ and $n_A \times (k+l)$ matrices. If enough computational resources are available, these can be parallelised. In the best case scenario, all $k+l$ matrix-vector products can be performed at the same time, making the cost of the matrix-matrix product equivalent to the cost of one iteration of CG plus communication between the processors. %, equivalent to $k+l$ matrix-vector products with $\bA$

When a randomised method is used to generate the preconditioner, an inner loop is performed as follows. Estimates of the Ritz values of the Hessian and the corresponding Ritz vectors are obtained with a randomised method (Algorithm~\ref{alg:revd}, \ref{alg:nystrom} or \ref{alg:ritzit}) and used to construct an LMP. Then the system \eqref{eq:Ax_eq_b} with the exact Hessian $\bold{A}$ is solved with PCG (Algorithm~\ref{alg:PCG}) using the LMP. The state is updated in the outer loop using the PCG solution.%$\delta \tilde{\bold{p}}$.}

\begin{table}
\resizebox{.85\linewidth}{!}{
\begin{minipage}{\textwidth}
 \begin{tabular}{|c| c| c|c|c|} 
 \hline
& Lanczos & REVD & Nystr\"{o}m  & REVD\_\textit{ritzit}   \\ 
 \hline
 Information source & Previous inner loop & Current inner loop & Current inner loop & Current inner loop\\
  \hline
\begin{tabular}{c}  Preconditioner for the \\  first inner loop  \end{tabular}  & No & Yes & Yes & Yes \\ 
  \hline
  \begin{tabular}{c}  $k$ dependency on the  \\  previous inner loop \end{tabular}  &  \begin{tabular}{c}  Bounded by the number \\ of CG iterations \end{tabular}  &Independent & Independent & Independent \\
   \hline
 \begin{tabular}{c}  Matrix-matrix products \\ with $\bA$ \end{tabular}  & None  &  \begin{tabular}{c} 2 products with \\ $n_A  \times (k+l)$ matrices  \end{tabular} & \begin{tabular}{c} 2 products with \\ $n_A  \times (k+l)$ matrices  \end{tabular} & \begin{tabular}{c}  1 product with \\ $n_A  \times (k+l)$ matrix \end{tabular}\\ % (parallelisable)
  \hline
 QR decomposition &  None & None & None  & $\breve{\bY} \in \mathbb{R}^{n_A  \times (k+l)}$  \\ 
\hline
 Orthogonalisation  &  None & $\bY \in \mathbb{R}^{n_A  \times (k+l)}$ &$\bY \in \mathbb{R}^{n_A  \times (k+l)}$ & $\bG \in \mathbb{R}^{n_A  \times (k+l)}$ \\ 
\hline
Cholesky factorization &  None & None &  $\bE^{(2)} \in \mathbb{R}^{ (k+l) \times (k+l) }$   & None  \\ 
\hline
Triangular solve & None & None & \begin{tabular}{c} $\bF \bC = \bE^{(1)}$ \\ for $\bF \in \mathbb{R}^{n_A  \times (k+l)}$\end{tabular}  & None \\
\hline
Deterministic EVD & $\bT_k \in \mathbb{R}^{k \times k}$ * & $\bar \bK \in \mathbb{R}^{ (k+l) \times (k+l) }$ & None & $\breve{\bK} \in \mathbb{R}^{ (k+l) \times (k+l) }$  \\
  \hline
  Deterministic SVD & None & None & $\bF \in \mathbb{R}^{n_A \times (k+l)}$ & None \\
  \hline
\end{tabular}
\end{minipage} }
\caption{A summary of the properties of the different methods of obtaining $k$ Ritz vectors and values to generate the preconditioner for a $n_A \times n_A$ matrix $\bA$ in the $i$th inner loop. Here $l$ is the oversampling parameter. * applies for CG with reorthogonalization. }
\label{table:all_formulations_summary}
\end{table}

\section{Numerical experiments}
We demonstrate our proposed preconditioning strategies using two models: a simple linear advection model to explore the spectra of the preconditioned Hessian and the nonlinear Lorenz 96 model \citep{Lorenz96} to explore the convergence of split preconditioned CG (PCG) . We perform identical twin experiments, where $\bold{x}^t = ((\bx^t_0)^T, \dots, (\bx^t_N)^T)^T$ denotes the reference trajectory. The observations and background state are generated by adding noise with covariance matrices $\bold{R}$ and $\bB$, respectively, to $\mathcal{H}_i (\bx^t_i)$ and $\bx_0$. We use direct observations, thus the observation operator $\mathcal{H}_i$ is linear. 

We use covariance matrices $\bR_i=\sigma_o^2 \bI_{q_i}$, where $q_i$ is the number of observations at time $t_i$, $\bQ_i = \sigma_q^2 \bC_q$, where $\bC_q$ is a Laplacian correlation matrix \citep{Johnson2005}, and $\bB=\sigma_b^2 \bC_b$, where $\bC_b$ is a second-order auto-regressive correlation matrix \citep{Daley91}.

%We consider the case when first level preconditioning is already applied, i.e. if we use no second level preconditioning then we solve \eqref{eq:forcing_1stlvl_prec}, and if we apply the second level preconditioning then we solve linear systems with preconditioned matrices $C_k^T A C_k$, where $C_k$ is an LMP constructed with $k$ vectors and $\bA$ is defined in \eqref{eq:forcing_1stlvl_prec}. 
We assume that first level preconditioning has already been applied (recall \eqref{eq:forcing_1stlvl_prec}). In data assimilation, using Ritz-LMP as formulated in \eqref{eq:Ritz_LMP} is impractical because of the matrix products with $\bA$ and we cannot use a simple formulation of Ritz-LMP when the Ritz values and vectors are obtained with the randomised methods. Hence, we use the spectral-LMP. However, as we mentioned in Section~\ref{sec:LMPs}, the spectral-LMP that is constructed with well converged Ritz values and vectors acts like Ritz-LMP. When we consider the second inner loop, we compare the spectral-LMPs with information from the randomised methods with the spectral-LMP constructed with information obtained with the Matlab function \textit{eigs} in the previous inner loop. \textit{eigs} returns a highly accurate estimate of a few largest or smallest eigenvalues and corresponding eigenvectors. We will use the term randomised LMP to refer to the spectral-LMPs that are constructed with information from the randomised methods, and deterministic LMP to refer to the spectral-LMP that is constructed with information from \textit{eigs}.

The computations are performed with Matlab R2017b. Linear systems are solved using the Matlab implementation of PCG (function \textit{pcg}), which was modified to allow split preconditioning to maintain the symmetric coefficient matrix at every loop.

\subsection{Advection model}
First, we consider the linear advection model: 
\begin{equation}\label{eq:advection}
\frac{\partial u(z,t)}{\partial t} + \frac{\partial u(z,t)}{\partial z} = 0,
\end{equation}
where $z \in [0,1]$ and  $t \in (0, T)$. An upwind numerical scheme is used to discretise \eqref{eq:advection} (see, e.g. Chapter 4 of \cite{Morton_Numsol}). To allow us to compute all the eigenvalues (described in the following section), we consider a small system with the linear advection model. The domain is divided into $n=40$ equally spaced grid points, with grid spacing $\Delta z = 1 / n$. We run the model for 51 time steps ($N=50$) with the time step size $\Delta t = 1 / N$, hence $\bA$ is a $2040 \times 2040$ matrix. The Courant number is $C = 0.8$ (the upwind scheme is stable with $C \in [0,1]$). The initial conditions are Gaussian $u(z,0) = 6  \exp \left( -\frac{(z-0.5)^2}{2\times 0.1^2} \right)$, and the boundary conditions are periodic $u(1,t) = u(0,t)$.
 
We set $\sigma_o = 0.05$, $\sigma_q = 0.05$ and $\sigma_b = 0.1$. $\bC_q$ and $\bC_b$ have length scales equal to $10 \Delta y$. Every 4th model variable is observed at every 5th time step, ensuring that there is an observation at the final time step (100 observations in total). Because the model and the observational operator $\mathcal{H}_i$ are linear, the cost function \eqref{eq:4D-var_error} is quadratic and its minimiser is found in the first loop of the incremental method. 

\subsubsection{Eigenvalues of the preconditioned matrix}
We apply the randomised LMPs in the first inner loop. Note that if the deterministic LMP is used, it is unclear how to precondition the first inner loop. We explore what effect the randomised LMPs have on the eigenvalues of $\bA$. The oversampling parameter $l$ is set to 5 and the randomised LMPs are constructed with $k=25$ vectors. 

\begin{figure}[h!]
\begin{subfigure}[b]{0.5\linewidth}
  \centering
 \includegraphics[width=\linewidth]{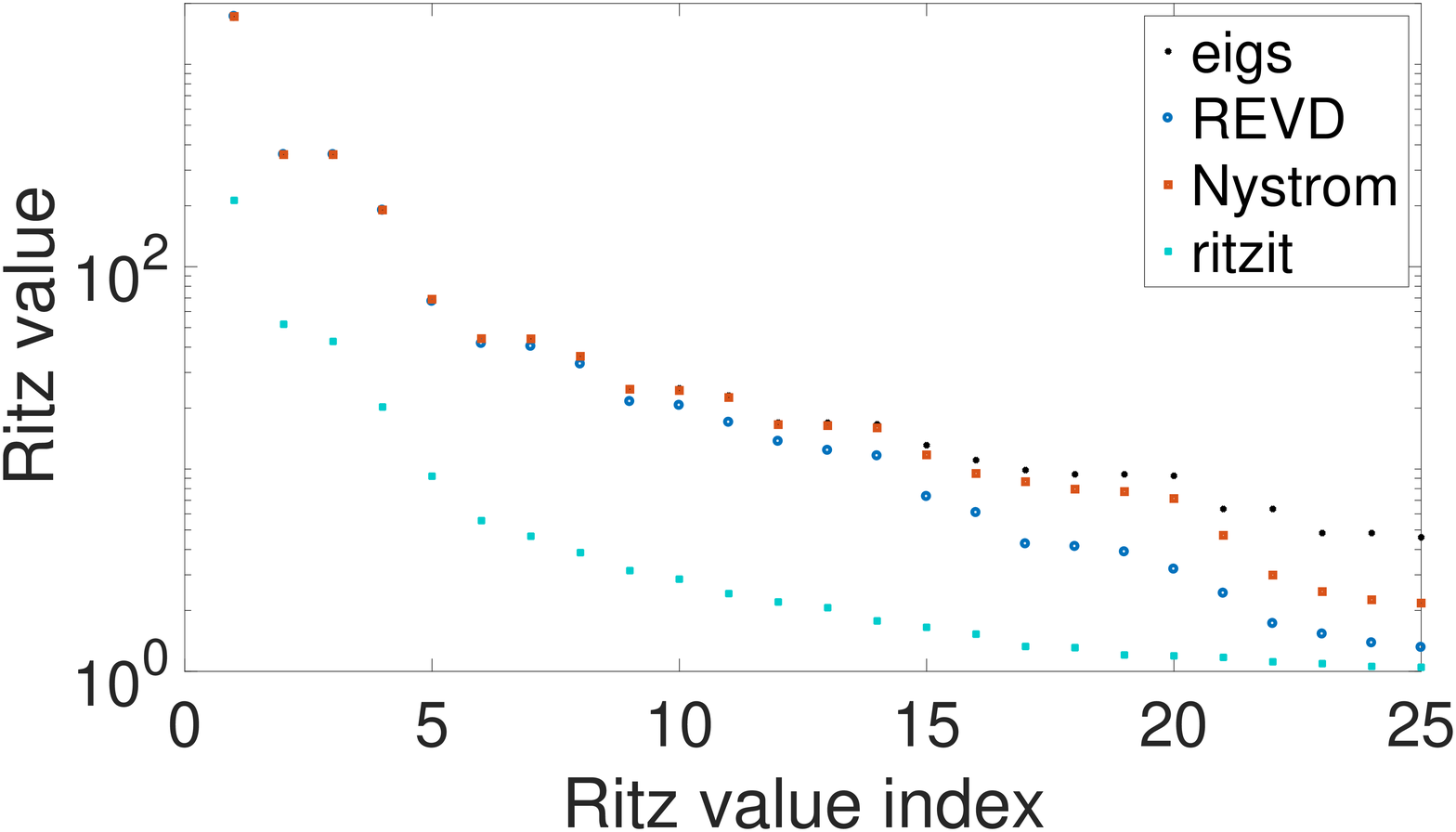}
     \caption{Ritz values} \label{fig:ritz-val-advection}
\end{subfigure}
\begin{subfigure}[b]{0.5\linewidth}
  \centering
 \includegraphics[width=\linewidth]{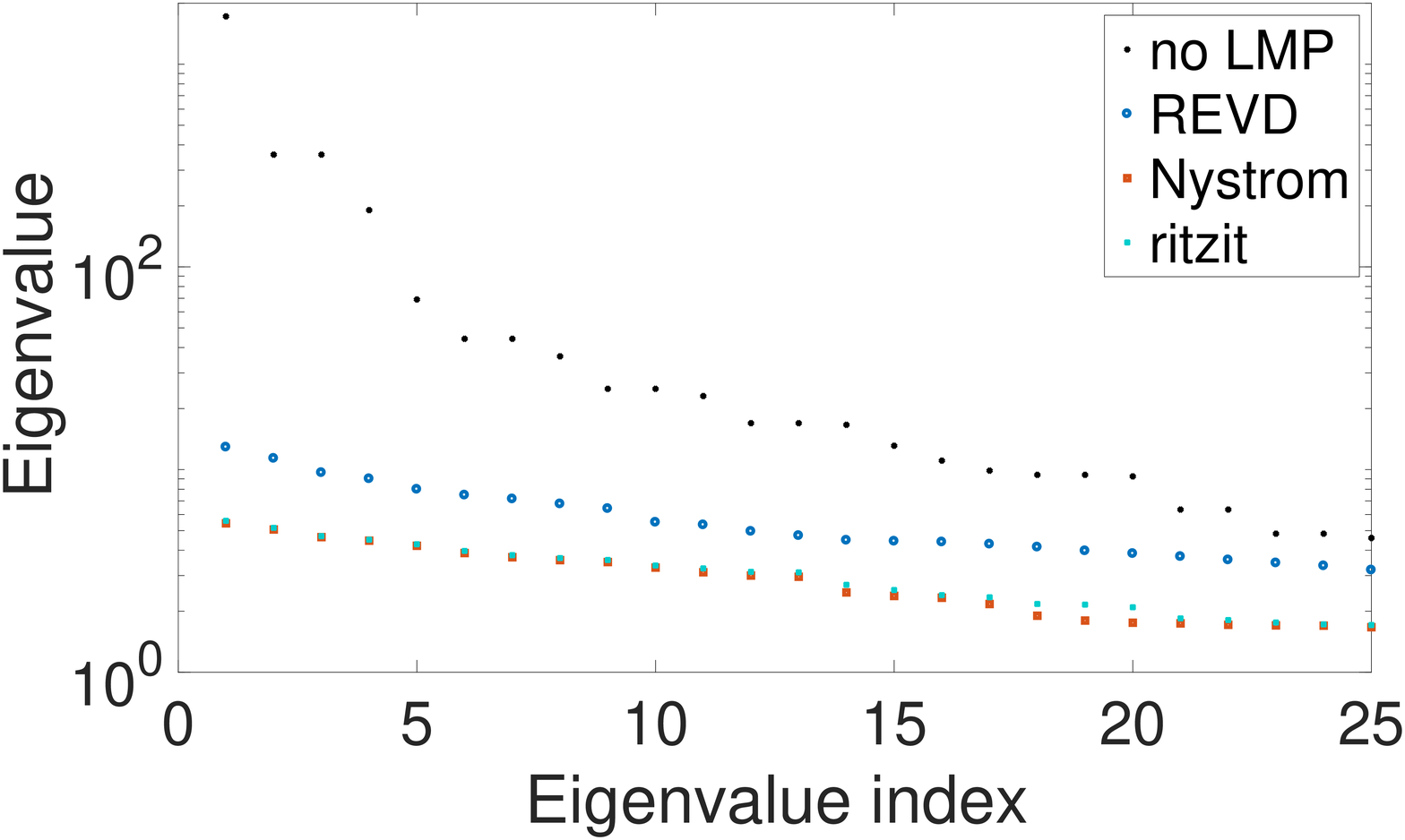}
   \caption{Preconditioned spectra, largest} \label{fig:preconditioned-eigs-advection-largest}
\end{subfigure}\\[1ex]
\begin{subfigure}[b]{0.5\linewidth}
  \centering
 \includegraphics[width=\linewidth]{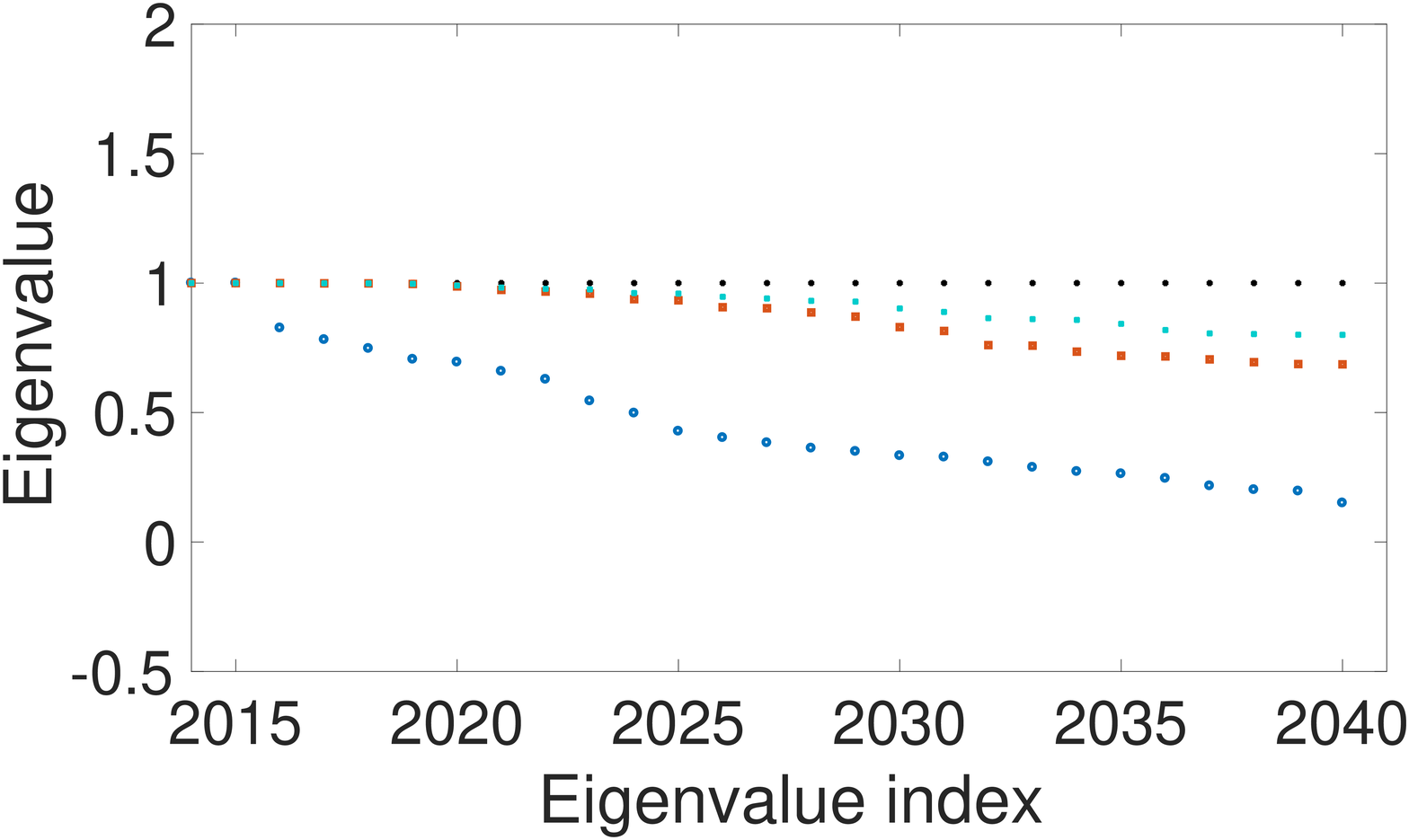}
     \caption{Preconditioned spectra, smallest} \label{fig:preconditioned-eigs-smallest-advection}
\end{subfigure}
\begin{subfigure}[b]{0.5\linewidth}
  \centering
 \includegraphics[width=\linewidth]{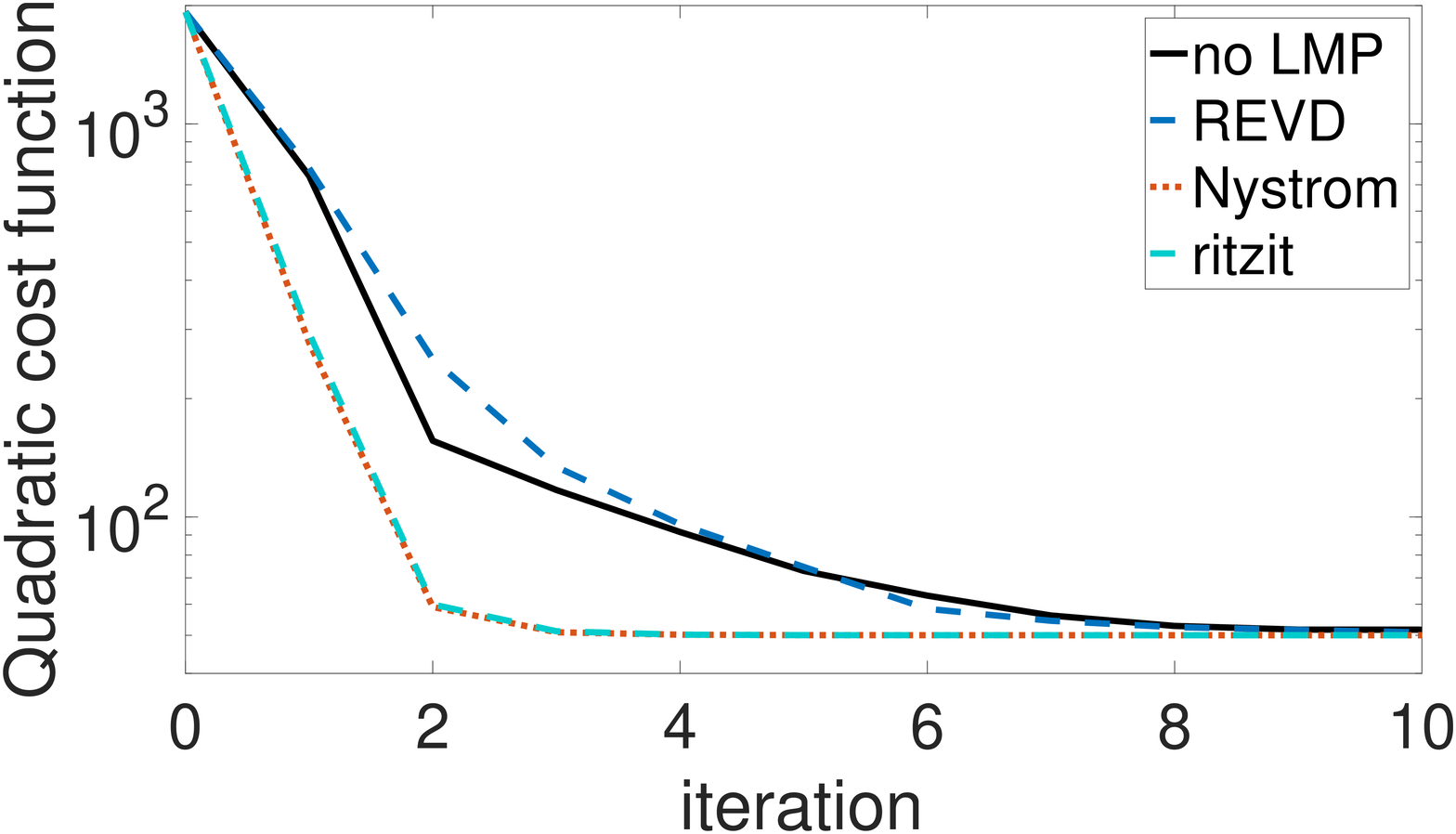}
     \caption{Quadratic cost function} \label{fig:qcf-vs-PCG_it-advection}
\end{subfigure}
\caption{Advection problem. (\subref{fig:ritz-val-advection}) 25 largest eigenvalues of $\bA$ (\textit{eigs}) and their estimates given by randomised methods, the largest eigenvalues and their estimates given by REVD and Nystr\"{o}m coincide. (\subref{fig:preconditioned-eigs-advection-largest}) largest eigenvalues of $\bA$ (no LMP, the same as \textit{eigs} in (\subref{fig:ritz-val-advection})) and $(\bC_{25}^{sp})^T \bA \bC_{25}^{sp}$, where $\bC_{25}^{sp}$ is constructed with Ritz values in (\subref{fig:ritz-val-advection}) and corresponding Ritz vectors. (\subref{fig:preconditioned-eigs-smallest-advection}) smallest eigenvalues of $\bA$ and $(\bC_{25}^{sp})^T \bA \bC_{25}^{sp}$. (\subref{fig:qcf-vs-PCG_it-advection}) quadratic cost function value versus PCG iteration when solving systems with $\bA$ and $(\bC_{25}^{sp})^T \bA \bC_{25}^{sp}$.} 
\label{fig:Ritz_values_advection_ps43}
\end{figure}

The Ritz values of $\bA$ given by the randomised methods are compared to those computed using \textit{eigs} (Figure~\ref{fig:ritz-val-advection}). The Nystr\"{o}m method produces a good approximation of the largest eigenvalues, REVD gives a slightly worse approximation, except for the largest five eigenvalues. The REVD\_\textit{ritzit} method underestimates the largest eigenvalues significantly. %The largest and smallest eigenvalues of $\bA$ and the preconditioned matrices are plotted in Figure~\ref{fig:preconditioned-eigs-advection-largest} and \ref{fig:preconditioned-eigs-smallest-advection}. 
The largest eigenvalues of the preconditioned matrices are smaller than the largest eigenvalue of $\bA$ (Figure~\ref{fig:preconditioned-eigs-advection-largest}). However, the smallest eigenvalues of the preconditioned matrices are less than one and hence applying the preconditioner expands the spectrum of $\bA$ at the lower boundary (Figure~\ref{fig:preconditioned-eigs-smallest-advection}), so that Theorem 3.4 of \cite{Gratton11}, which considers the non-expansiveness of the spectrum of the Hessian after preconditioning with an LMP, does not hold. This happens because the formulation of the spectral-LMP is derived assuming that the eigenvalues and eigenvectors are exact, while the randomized methods provide only approximations. Note that even though REVD\_\textit{ritzit} gives the worst approximations of the largest eigenvalues of the Hessian, using the randomised LMP with information from REVD\_\textit{ritzit} reduces the largest eigenvalues of the preconditioned matrix the most and the smallest eigenvalues are close to one. Using the randomised LMP with estimates from Nystr\"{o}m gives similar results. Hence, the condition number of the preconditioned matrix is lower when the preconditioners are constructed with REVD\_\textit{ritzit} or Nystr\"{o}m compared to REVD.

The values of the quadratic cost function at the first ten iterations of PCG are shown in Figure~\ref{fig:qcf-vs-PCG_it-advection}. Using the randomised LMP that is constructed with information from REVD is detrimental to the PCG convergence compared to using no preconditioning. Using information from the Nystr\"{o}m and REVD\_\textit{ritzit} methods results in similar PCG convergence and low values of the quadratic cost function are reached in fewer iterations than without preconditioning. The PCG convergence may be explained by the favourable distribution of the eigenvalues after preconditioning using Nystr\"{o}m and REVD\_\textit{ritzit}, and the smaller than one eigenvalues when using REVD. These results, however, do not necessarily generalize to an operational setting as this system is well conditioned while operational settings are not. This will be investigated further in the next section.

\subsection{Lorenz 96 model}

We next use the Lorenz 96 model to examine what effect the randomised LMPs have on PCG performance. In the Lorenz 96 model the evolution of the $n$ components $X^j, \ j \in \{1,2,\dots,n\}$ of $\bx_i$ is governed by a set of $n$ coupled ODEs:
\begin{equation}\label{eq:lorenz96}
\frac{dX^j}{dt} = -X^{j-2} X^{j-1} + X^{j-1} X^{j+1} - X^j + F,
\end{equation}
where $X^{-1} = X^{n-1}, X^0 = X^n$ and $X^{n+1} = X^1$ and $F=8$. The equations are integrated using a fourth order Runge-Kutta scheme \citep{Butcher87}. We set $n=80$ and $N=150$ (the size of $\bA$ is $12080 \times 12080$) and observe every 10th model variable at every 10th time step (120 observations in total), ensuring that there are observations at the final time step. The grid point distance is $\Delta X = 1/n$ and the time step is set to $\Delta t = 2.5 \times 10^{-2}$. 

For the covariance matrices we use $\sigma_o =0.15 $ and $\sigma_b = 0.2$. $\bC_b$ has length scale equal to $2 \Delta X$. Two setups are used for the model error covariance matrix:
\begin{itemize}
\item $\sigma_q = 0.1$ and $\bC_q$ has length scale $L_q = 2 \Delta X$ (the same as $\bC_b$); %ps37
\item $\sigma_q = 0.05$ and $\bC_q$ has length scale $L_q = 0.25 \Delta X$. %ps64
\end{itemize}
In our numerical experiments, the preconditioners have very similar effect using both setups. Hence, we present plots for the case $\sigma_q = 0.1$ and $L_q = 2 \Delta X$ in the following sections, except Figure~\ref{fig:qcf_previous_ritzit_means_ft10_fx10_ps37-64_k-5-10-15}.

The first outer loop is performed and no second level preconditioning is used in the first inner loop, where PCG is run for 100 iterations or until the relative residual norm reaches $10^{-6}$. In the following sections, we use randomised and deterministic LMPs in the second inner loop. PCG has the same stopping criteria as in the first inner loop.

\subsubsection{Minimising the inner loop cost function}\label{sec:numerics_lorenz_cf}
\begin{figure}[h!]
\begin{subfigure}[b]{0.5\linewidth}
  \centering
 \includegraphics[width=\linewidth]{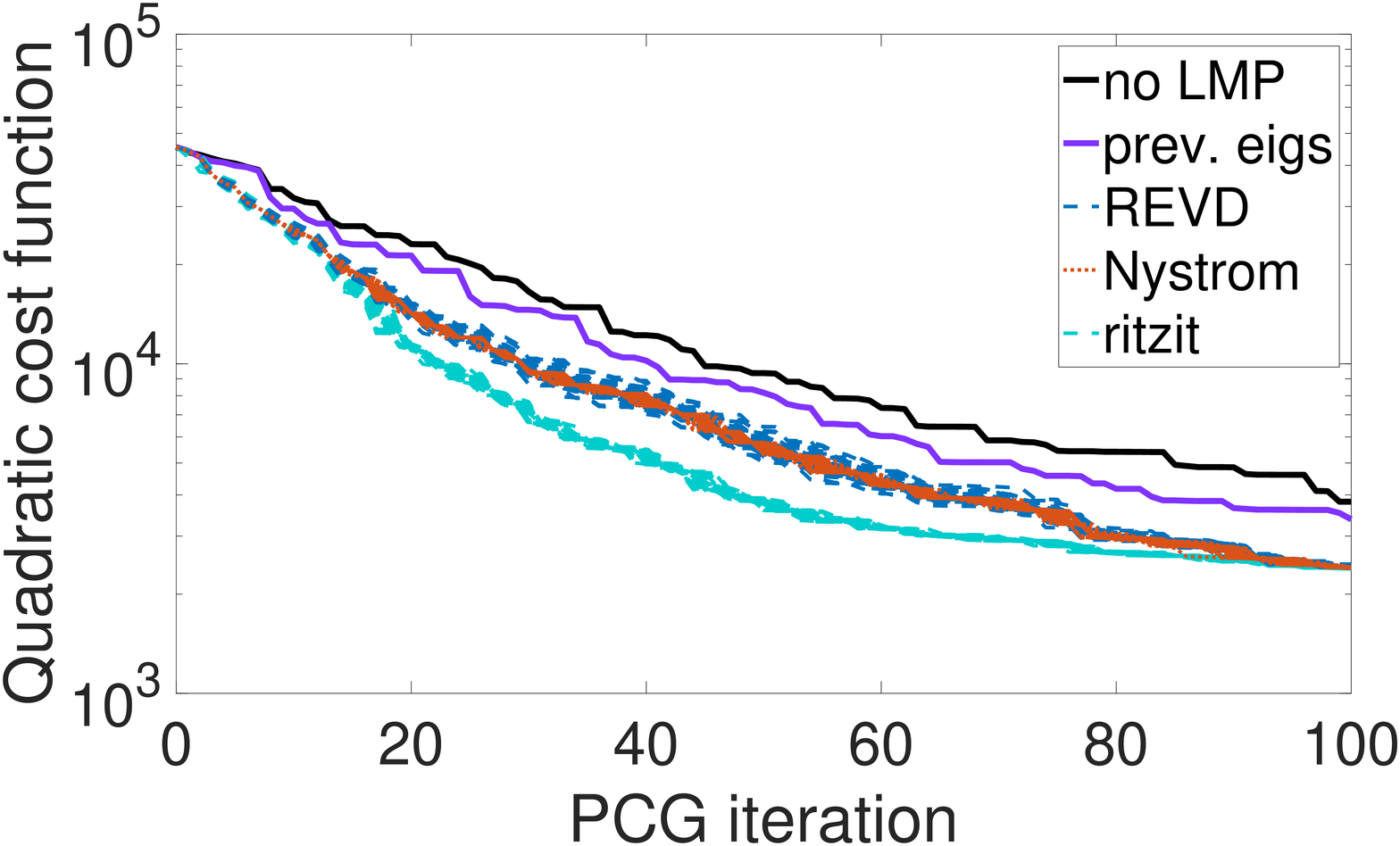}
     \caption{$k=5$, all runs} \label{fig:qcf_all_runns_ft10_fx10_ps37_k5}
\end{subfigure}
\begin{subfigure}[b]{0.5\linewidth}
  \centering
 \includegraphics[width=\linewidth]{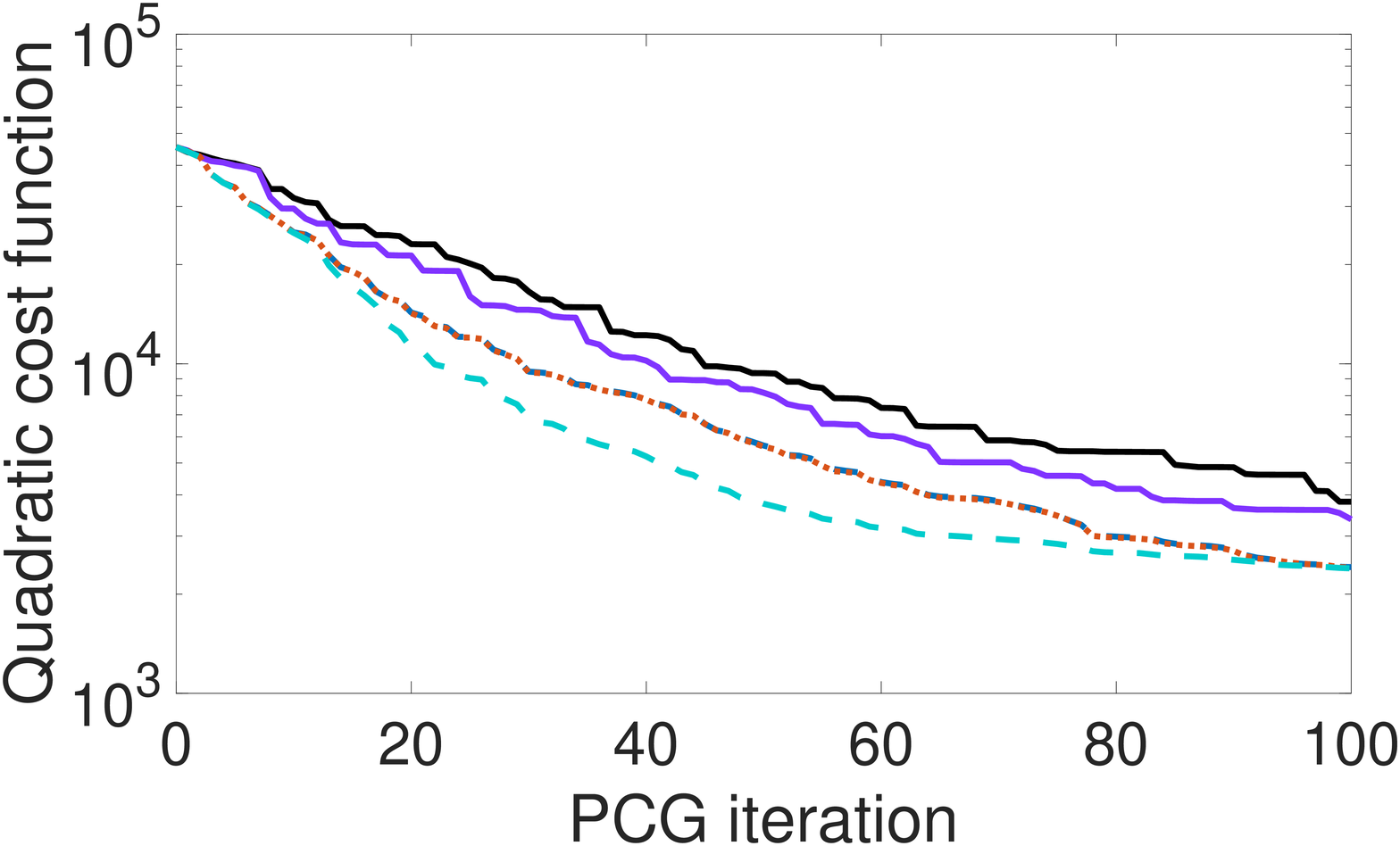}
     \caption{$k=5$} \label{fig:qcf_means_ft10_fx10_ps37_k5}
\end{subfigure}\\[1ex]
\begin{subfigure}[b]{0.5\linewidth}
  \centering
 \includegraphics[width=\linewidth]{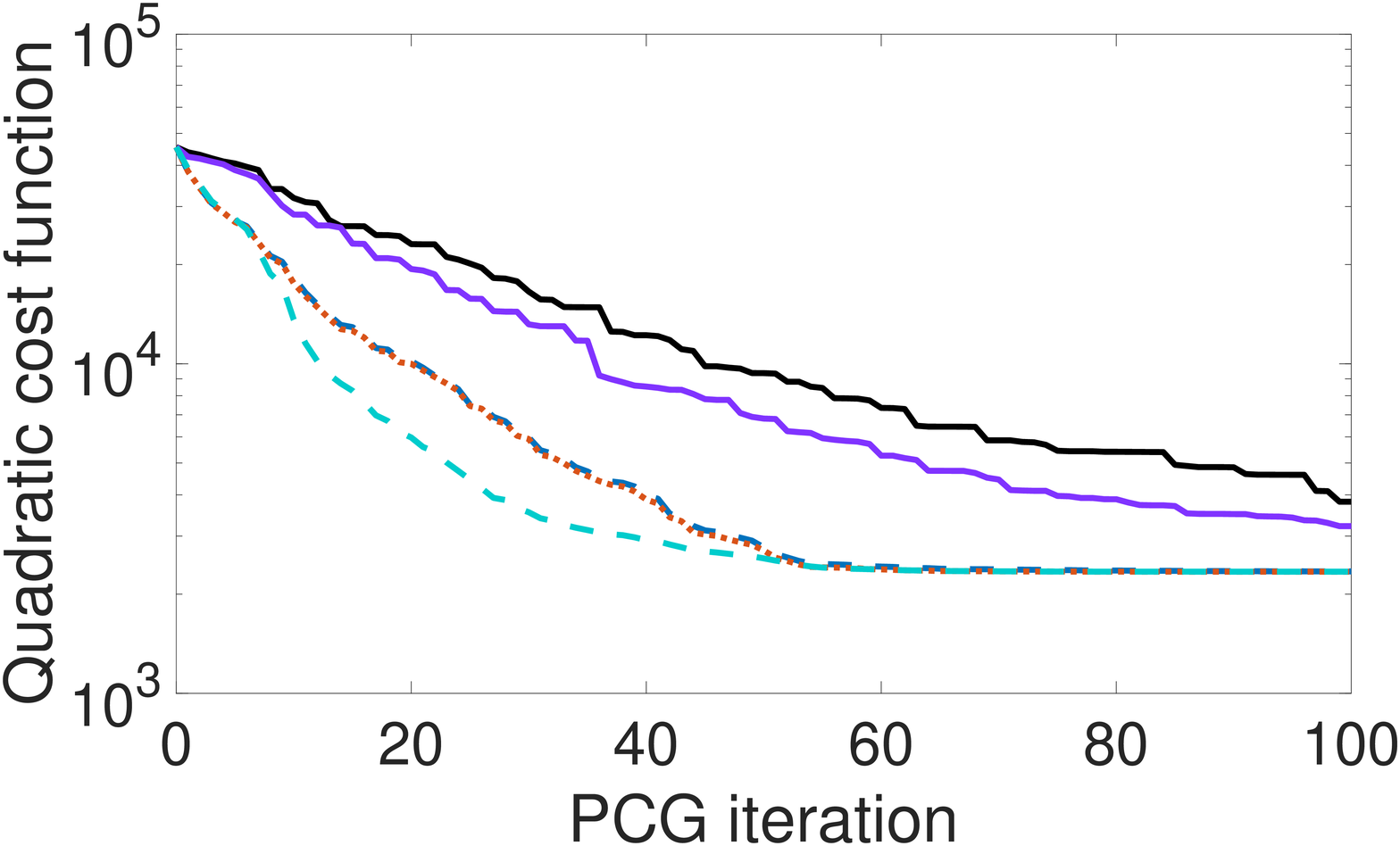}
   \caption{$k=10$} \label{fig:qcf_means_ft10_fx10_ps37_k10}
\end{subfigure}
\begin{subfigure}[b]{0.5\linewidth}
  \centering
 \includegraphics[width=\linewidth]{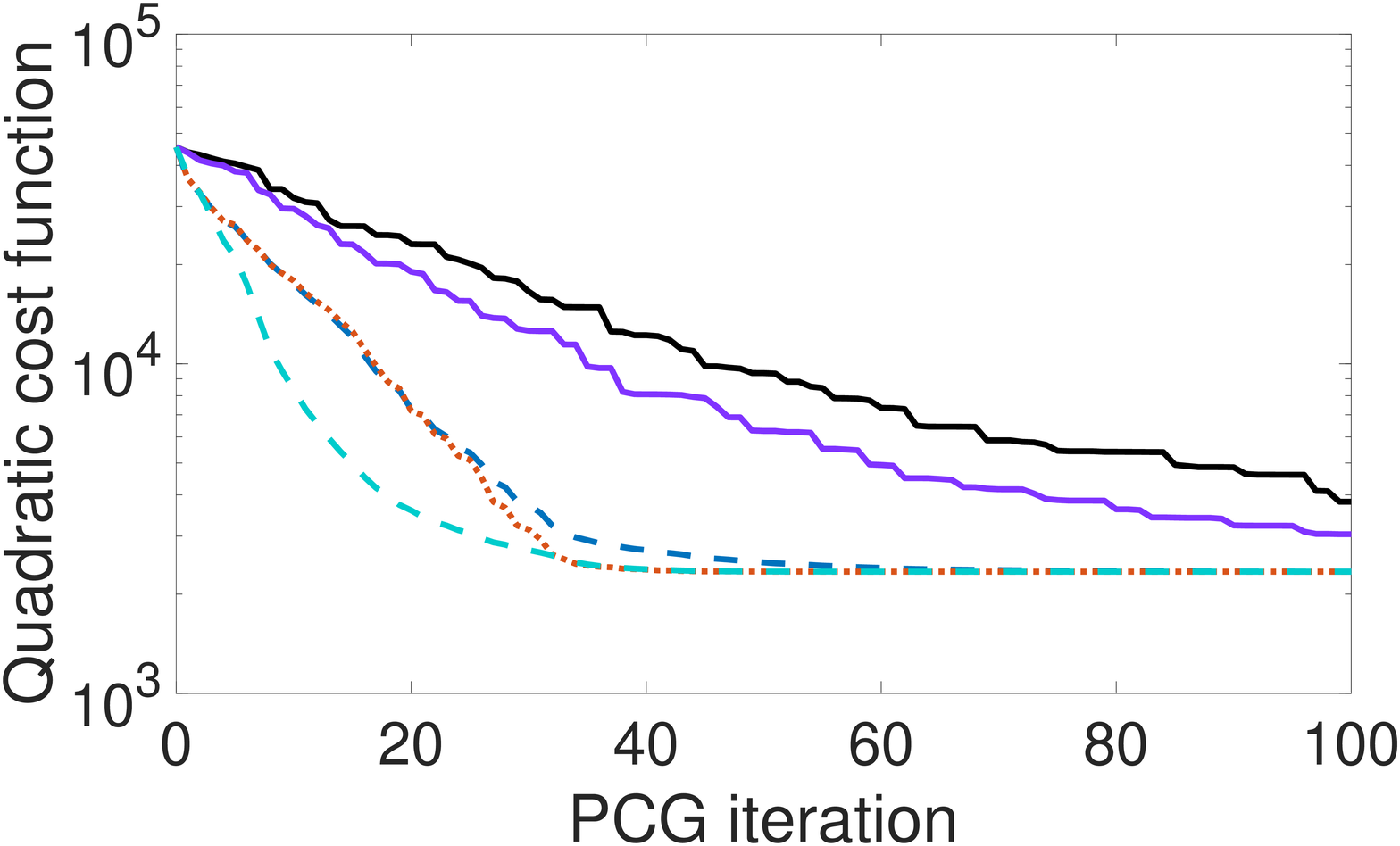}
     \caption{$k=15$} \label{fig:qcf_means_ft10_fx10_ps37_k15}
\end{subfigure}
\caption{A comparison of the value of the quadratic cost function at every PCG iteration when the spectral-LMP is constructed with $k \in \{5, 10, 15\}$ Ritz values and vectors obtained with the randomised methods in the current inner loop, and function \textit{eigs} in the previous inner loop. We also show no second level preconditioning (no LMP), which is the same in all four panels. For the randomised methods, (\subref{fig:qcf_all_runns_ft10_fx10_ps37_k5}) shows 50 experiments for $k=5$ and the rest display means over 50 experiments. Here $\sigma_q = 0.1$ and $L_q = 2 \Delta X$.}
\label{fig:qcf_means_ft10_fx10_ps37_k-5-10-15}
\end{figure}

In Figure~\ref{fig:qcf_means_ft10_fx10_ps37_k-5-10-15}, we compare the performance of the randomised LMPs %when we use Ritz values and vectors from the current inner loop 
with the deterministic LMP. %constructed with Ritz values and vectors that are obtained in the previous loop with \textit{eigs} function. 
We also consider the effect of varying $k$, the number of vectors used to construct the preconditioner. We set the oversampling parameter to $l=5$. %In Figure~\ref{fig:qcf_means_ft10_fx10_ps37_k-5-10-15}, we compare the value of the quadratic cost function at every PCG iteration when there is no second level preconditioning, the randomised LMPs are constructed with Ritz values and vectors obtained with REVD, Nystr\"{o}m and REVD\_$ritzit$ methods ($l=5$), and when information from the previous inner loop is used in the deterministic LMP. 
Because results from the randomized methods depend on the random matrix used, we perform 50 experiments with different realizations for the random matrix. We find that the different realizations lead to very similar results (see Figure~\ref{fig:qcf_all_runns_ft10_fx10_ps37_k5}).
%We construct fifty preconditioners for every randomised method by using fifty different random Gaussian matrices and solve the systems with these preconditioners. The results for every method are similar with all fifty random matrices.
 %and present the cost function value for $k=5$. The average cost function value at every iteration is displayed for all $k$ values (5, 10 and 15). 

Independently of the $k$ value, there is an advantage in using the second level preconditioning. The reduction in the value of the quadratic cost function is faster using the randomised LMPs compared to the deterministic LMPs, with REVD\_\textit{ritzit} performing the best after the first few iterations. The more information we use in the preconditioner (i.e. the higher $k$ value), the faster REVD\_\textit{ritzit} shows better results than the other methods. The performances of the REVD and Nystr\"{o}m methods are similar. Note that as $k$ increases, the storage (see Table~\ref{table:all_formulations_summary}) and work per PCG iteration increase. Examination of the Ritz values given by the randomised methods shows that REVD\_\textit{ritzit} gives the worse estimate of the largest eigenvalues, as was the case when using the advection model. We calculated the smallest eigenvalue of the preconditioned matrix $(\bC_{5}^{sp})^T \bA \bC_{5}^{sp}$ using \textit{eigs}. When $\bC_{5}^{sp}$ is constructed using REVD\_\textit{ritzit} or Nystr\"{o}m the smallest eigenvalue of $(\bC_{5}^{sp})^T \bA \bC_{5}^{sp}$ is equal to one, whereas using REVD it is approximately $0.94$. This may explain why the preconditioner constructed using REVD may not perform as well as other randomised preconditioners, but it is not entirely clear why the preconditioner that uses REVD\_\textit{ritzit} shows the best performance.

\begin{figure}[h]
\begin{subfigure}[b]{0.5\linewidth}
  \centering
 \includegraphics[width=\linewidth]{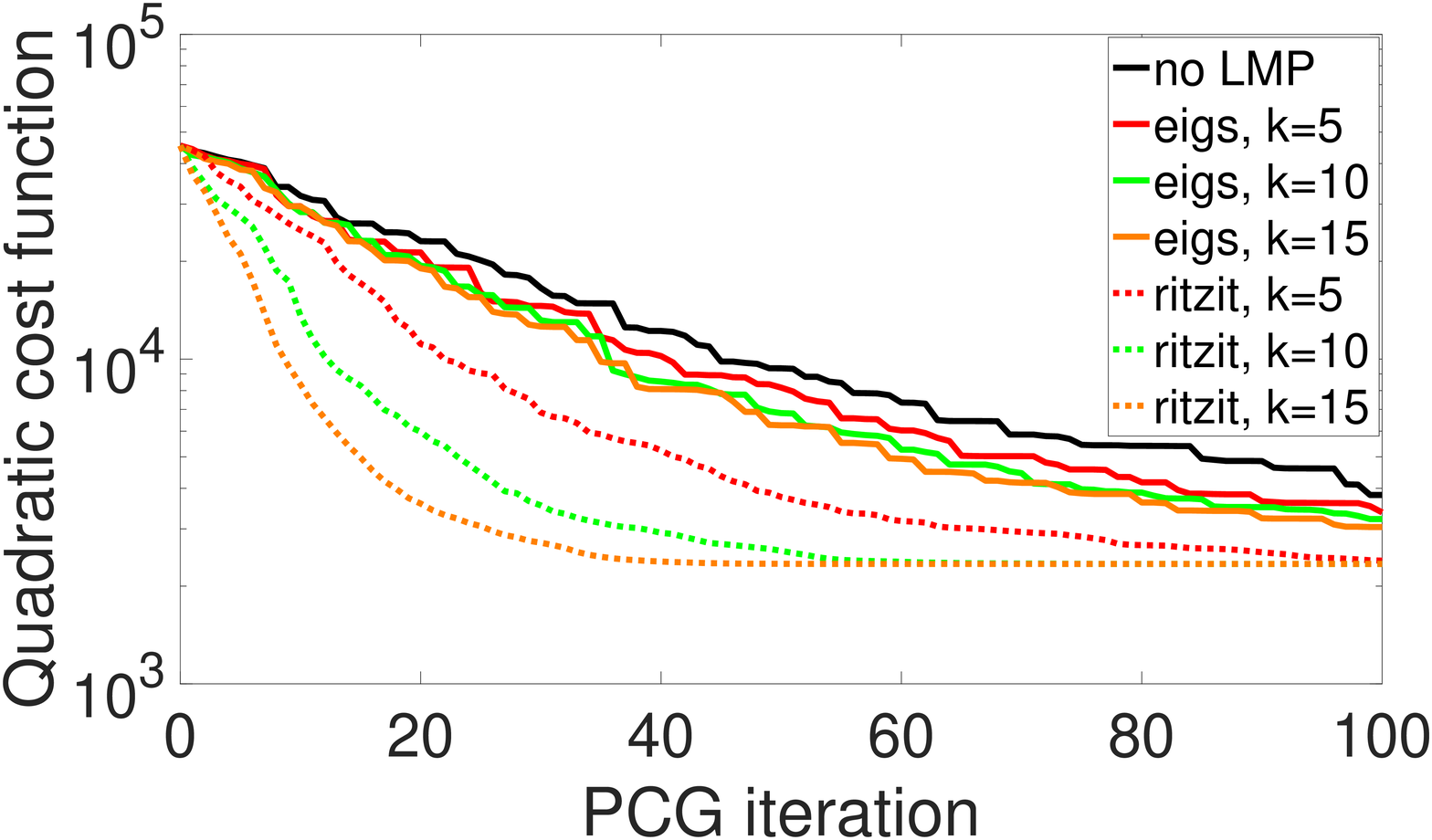}
     \caption{$\sigma_q = 0.1$ and $L_q = 2 \Delta X$} \label{fig:qcf_previous_ritzit_means_ft10_fx10_ps37_k-5-10-15}
\end{subfigure}
\begin{subfigure}[b]{0.5\linewidth}
  \centering
 \includegraphics[width=\linewidth]{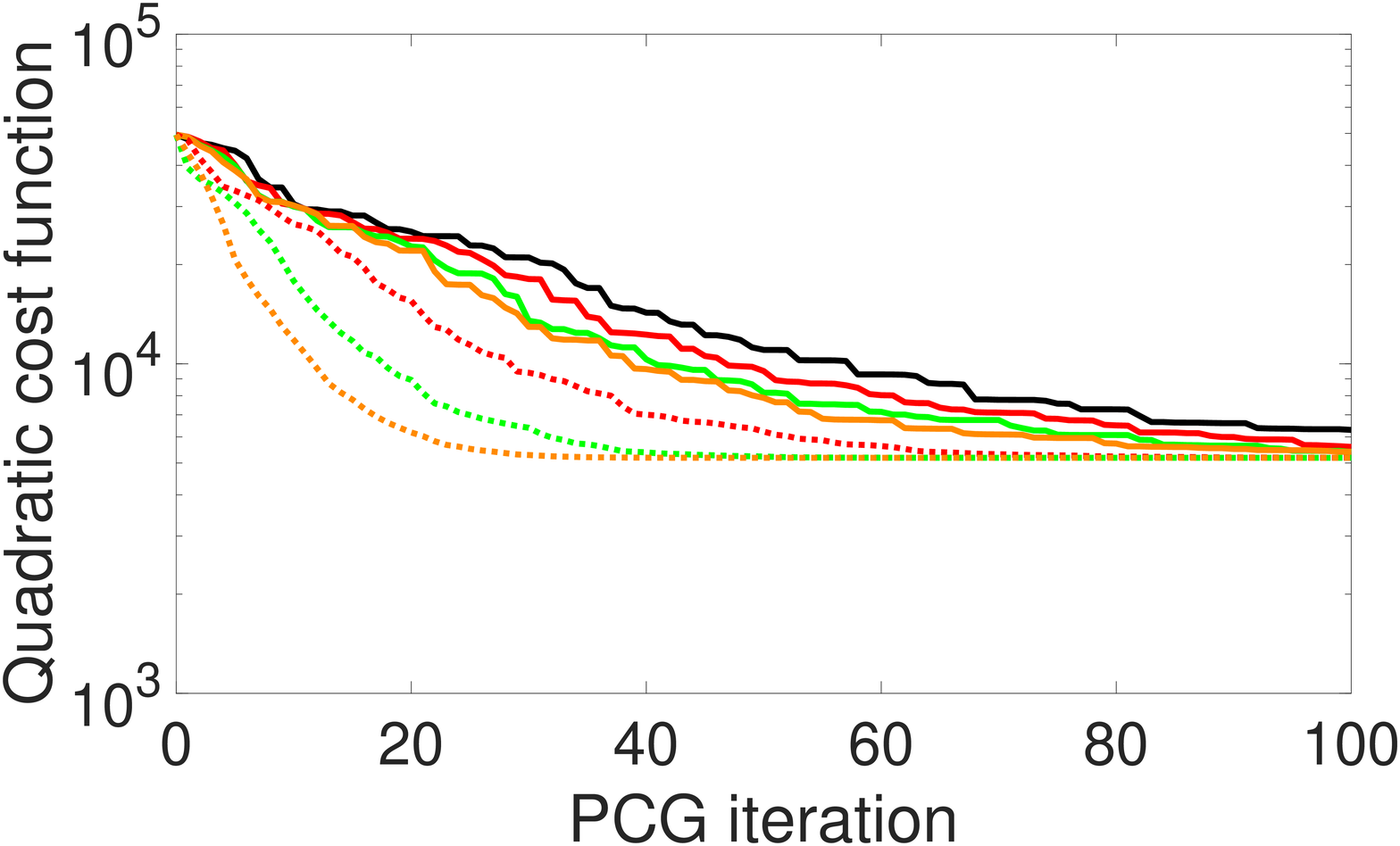}
     \caption{$\sigma_q = 0.05$ and $L_q = 0.25 \Delta X$} \label{fig:qcf_previous_ritzit_means_ft10_fx10_ps64_k-5-10-15}
\end{subfigure}\\[1ex]
\caption{A comparison of the values of the quadratic cost function at every PCG iteration when using deterministic LMP with information from the previous loop (eigs) and the randomised LMP with information from REVD\_\textit{ritzit} for different $k$ values (5, 10 and 15). No second level preconditioning is also shown (case (\subref{fig:qcf_previous_ritzit_means_ft10_fx10_ps37_k-5-10-15}) is the same as in Figure~\ref{fig:qcf_means_ft10_fx10_ps37_k-5-10-15}). In cases (\subref{fig:qcf_previous_ritzit_means_ft10_fx10_ps37_k-5-10-15}) and (\subref{fig:qcf_previous_ritzit_means_ft10_fx10_ps64_k-5-10-15}) the model error covariance matrices are constructed using parameters $\sigma_q $ and $L_q$.}
\label{fig:qcf_previous_ritzit_means_ft10_fx10_ps37-64_k-5-10-15}
\end{figure}

The PCG convergence when using the deterministic LMP and the randomised LMP with information from REVD\_\textit{ritzit} with different $k$ values is compared in Figure~\ref{fig:qcf_previous_ritzit_means_ft10_fx10_ps37-64_k-5-10-15} for both setups of the model error covariance matrix. %The deterministic LMP gives a small improvement if any when $k$ value is increased, especially in the first iterations of PCG. 
For the deterministic LMP, varying $k$ has little effect, especially in the first iterations of PCG. However, for REVD\_\textit{ritzit}, increasing $k$ results in a greater decrease of the cost function in the first iterations of PCG. Also, at any iteration of PCG we obtain a lower value of the quadratic cost function using the randomised LMP with $k=5$ compared to the deterministic LMP with $k=15$, which uses exact eigenpair information from the Hessian of the previous loop.

%from Figure~\ref{fig:qcf_means_ft10_fx10_ps37_k-5-10-15} plotted together 

\subsubsection{Effect of the observation network}

To understand the sensitivities of the results from the different LMPs to the observation network, %randomised LMPs and see how their effect on PCG changes with the observation network, 
%We explore if the effect of the spectral-LMPs on PCG depend on the observation network. 
we consider a system with the same parameters as in the previous section, where we had 120 observations, but we now observe
\begin{itemize}
\item every 5th model variable at every 5th time step (480 observations in total);
\item every 2nd variable at every 2nd time step (3000 observations in total). 
\end{itemize}
The oversampling parameter is again set to $l=5$ and we set $k=5$ and $k=15$ for both observation networks. Since the number of observations is equal to the number of eigenvalues that are larger than one and there are more observations than in the previous section, there are more eigenvalues that are larger than one after the first level preconditioning. %We perform the first inner loop as before and consider the second inner loop. 
Because all 50 experiments with different Gaussian matrices in the previous section were close to the mean, we perform 10 experiments for each randomised method, solve the systems and report the means of the quadratic cost function.

\begin{figure}[h!]
\begin{subfigure}[b]{0.5\linewidth}
  \centering
 \includegraphics[width=\linewidth]{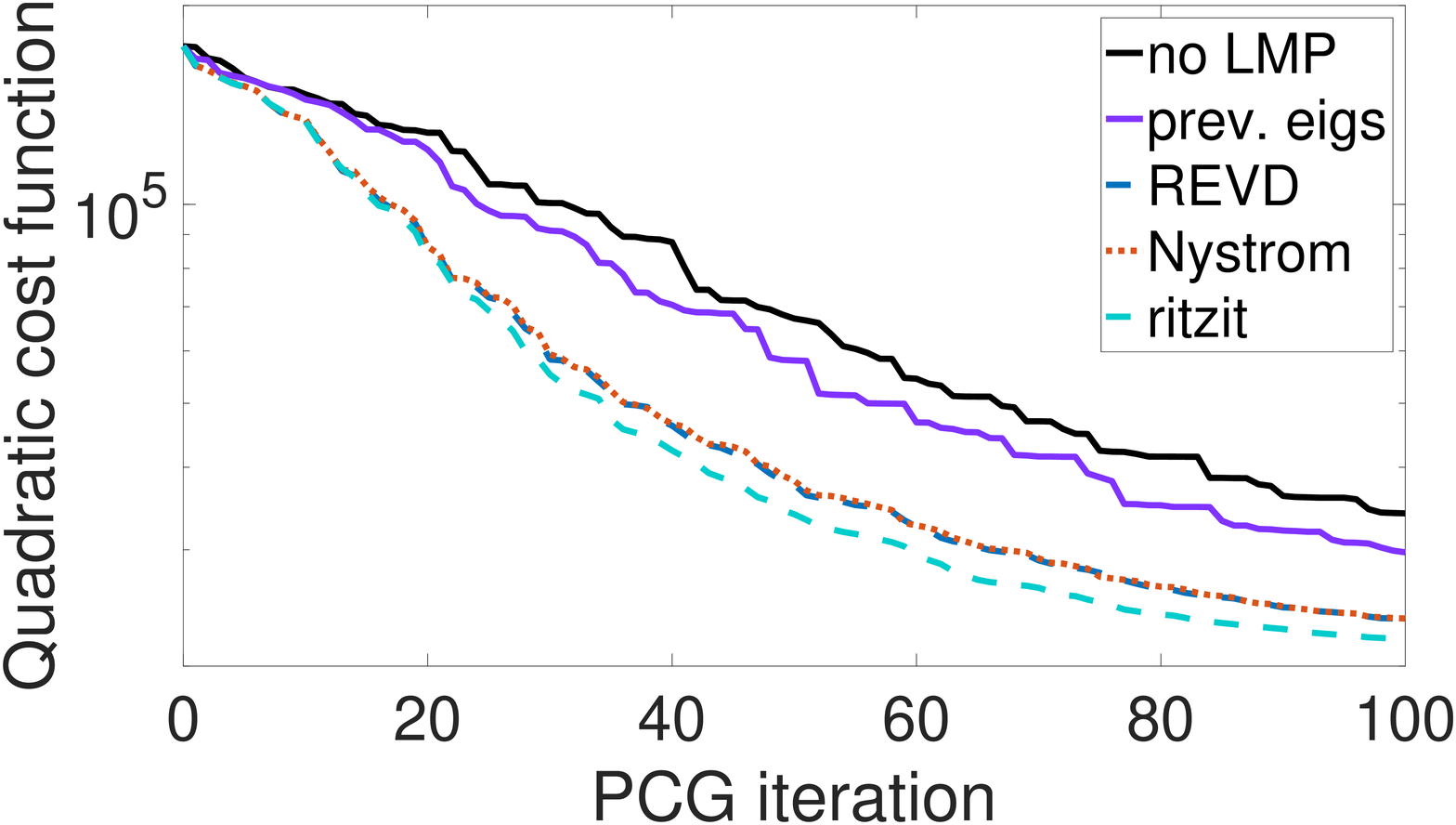}
     \caption{$k=5$, $q=480$} \label{fig:qcf_means_ft5_fx5_ps37_k5}
\end{subfigure}
\begin{subfigure}[b]{0.5\linewidth}
  \centering
 \includegraphics[width=\linewidth]{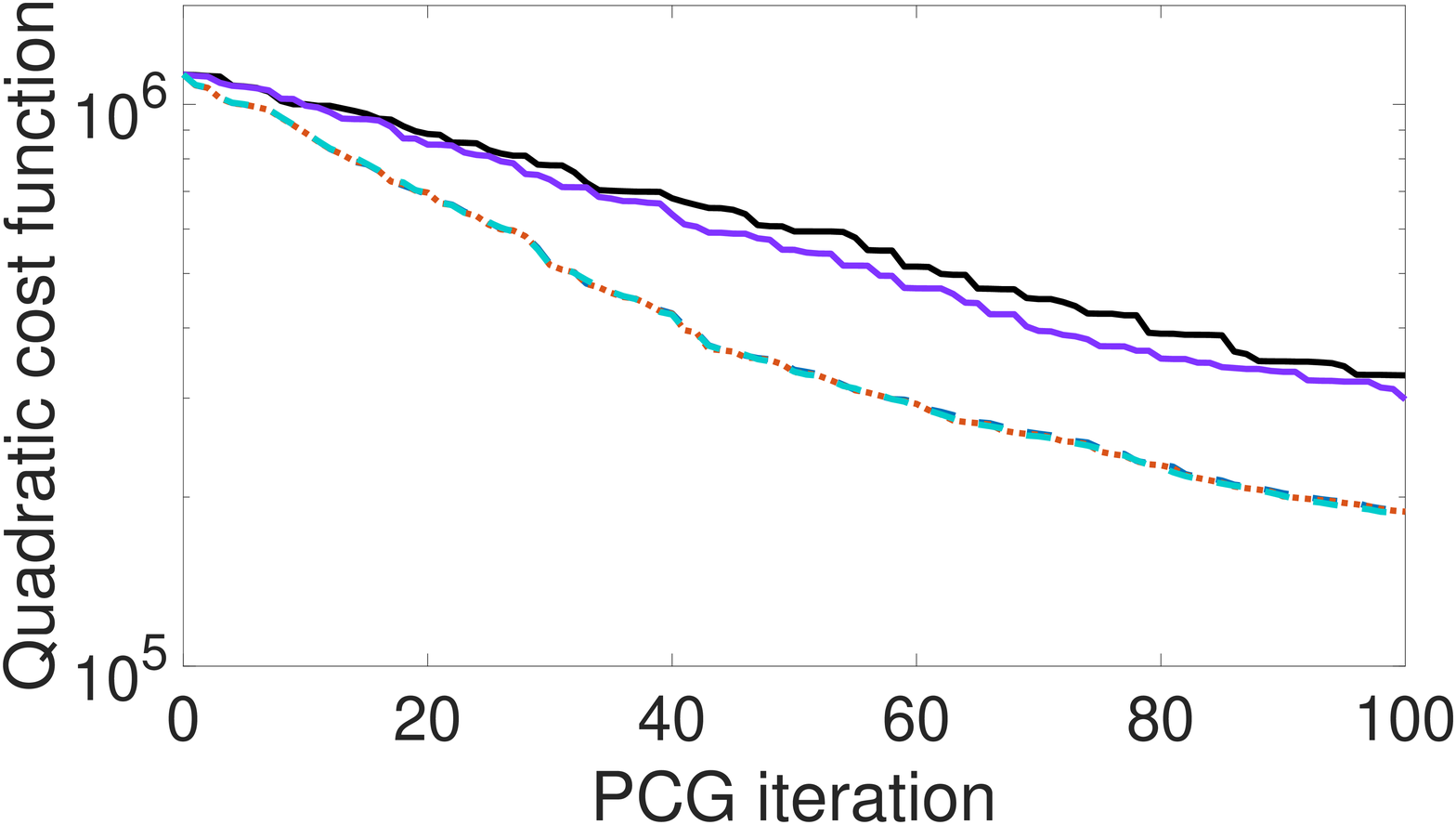}
     \caption{$k=5$, $q=3000$} \label{fig:qcf_means_ft2_fx2_ps37_k5}
\end{subfigure}\\[1ex]
\begin{subfigure}[b]{0.5\linewidth}
  \centering
 \includegraphics[width=\linewidth]{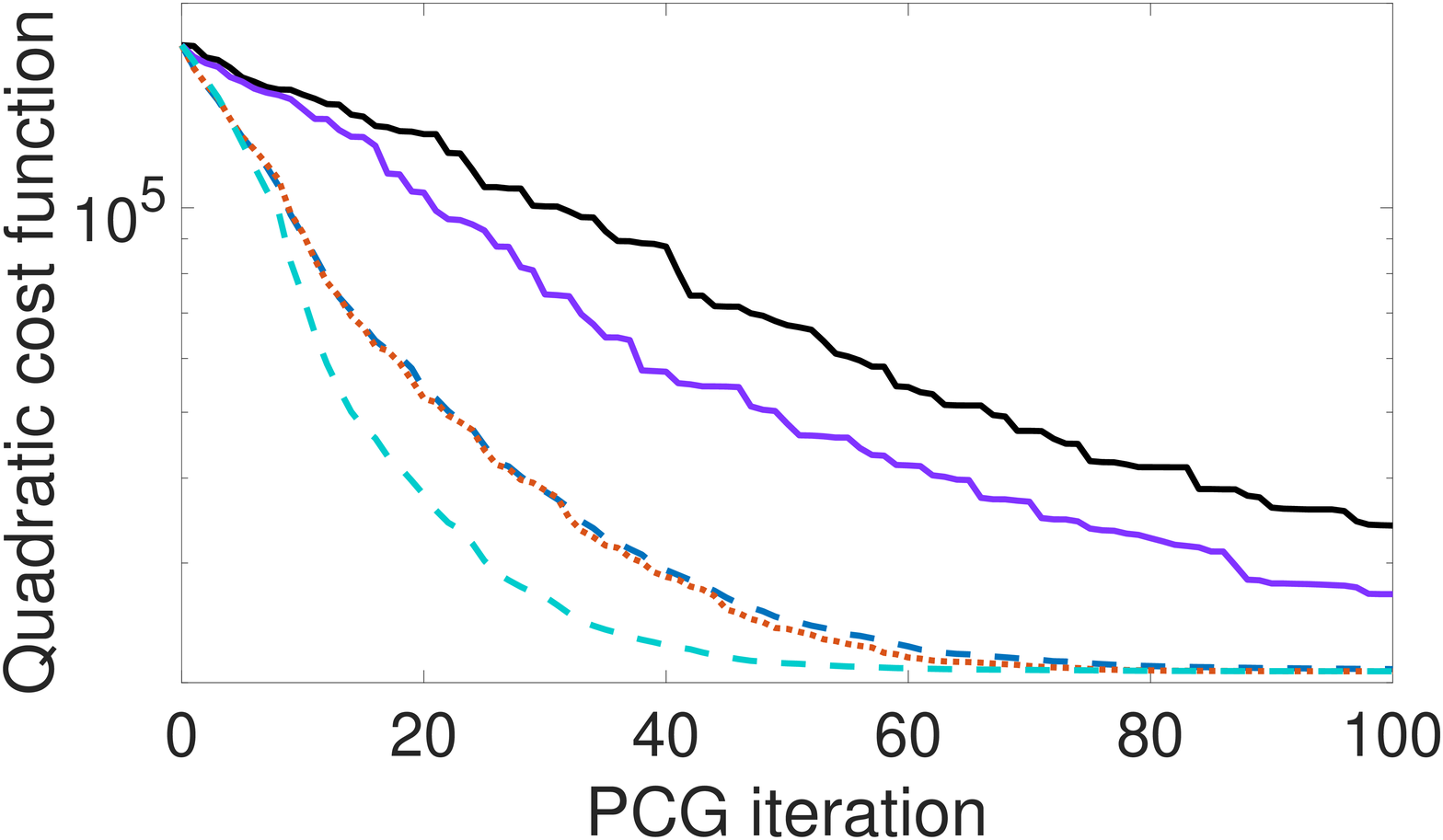}
     \caption{$k=15$, $q=480$} \label{fig:qcf_means_ft5_fx5_ps37_k15}
\end{subfigure}
\begin{subfigure}[b]{0.5\linewidth}
  \centering
 \includegraphics[width=\linewidth]{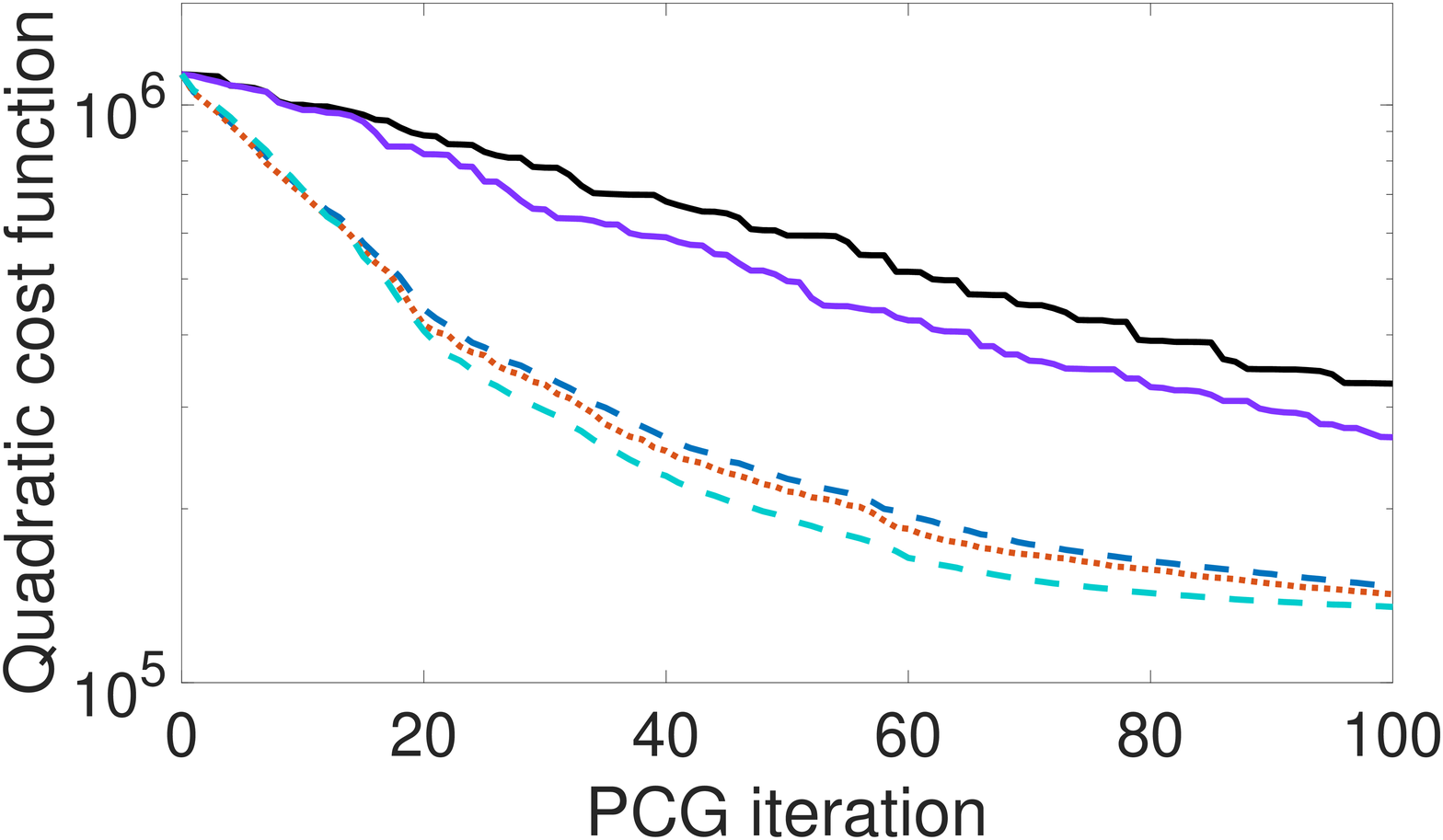}
     \caption{$k=15$, $q=3000$} \label{fig:qcf_means_ft2_fx2_ps37_k15}
\end{subfigure}
\caption{As in Figure~\ref{fig:qcf_means_ft10_fx10_ps37_k-5-10-15}, but for two systems with $q$ observations and 10 experiments are done for each randomised method and the mean values plotted.}
\label{fig:qcf_means_ft2_fx2_fx5_ft5_ps37_k-5-10-15}
\end{figure}

The results are presented in Figure~\ref{fig:qcf_means_ft2_fx2_fx5_ft5_ps37_k-5-10-15}. %As for the system with fewer observations in Figure~\ref{fig:qcf_previous_ritzit_means_ft10_fx10_ps37_k-5-10-15}, 
Again, the randomised LMPs perform better than the deterministic LMP. However, if the preconditioner is constructed with a small amount of information about the system ($k=5$ for both systems and $k=15$ for the system with 3000 observations), then there is little difference in the performance of different randomised LMPs. 
Also, when the number of observations is increased, more iterations of PCG are needed to get any improvement in the minimisation of the quadratic cost function when using the deterministic LMP over using no second level preconditioning. 

When comparing the randomised and deterministic LMPs with different values of $k$ for these systems, we obtain similar results to those in Figure~\ref{fig:qcf_previous_ritzit_means_ft10_fx10_ps37_k-5-10-15}, i.e. it is more advantageous to use the randomised LMP constructed with $k=5$ than using the deterministic LMP constructed with $k=15$.

\subsubsection{Effect of oversampling}
\begin{figure}[h!]
\begin{subfigure}[b]{0.5\linewidth}
  \centering
 \includegraphics[width=\linewidth]{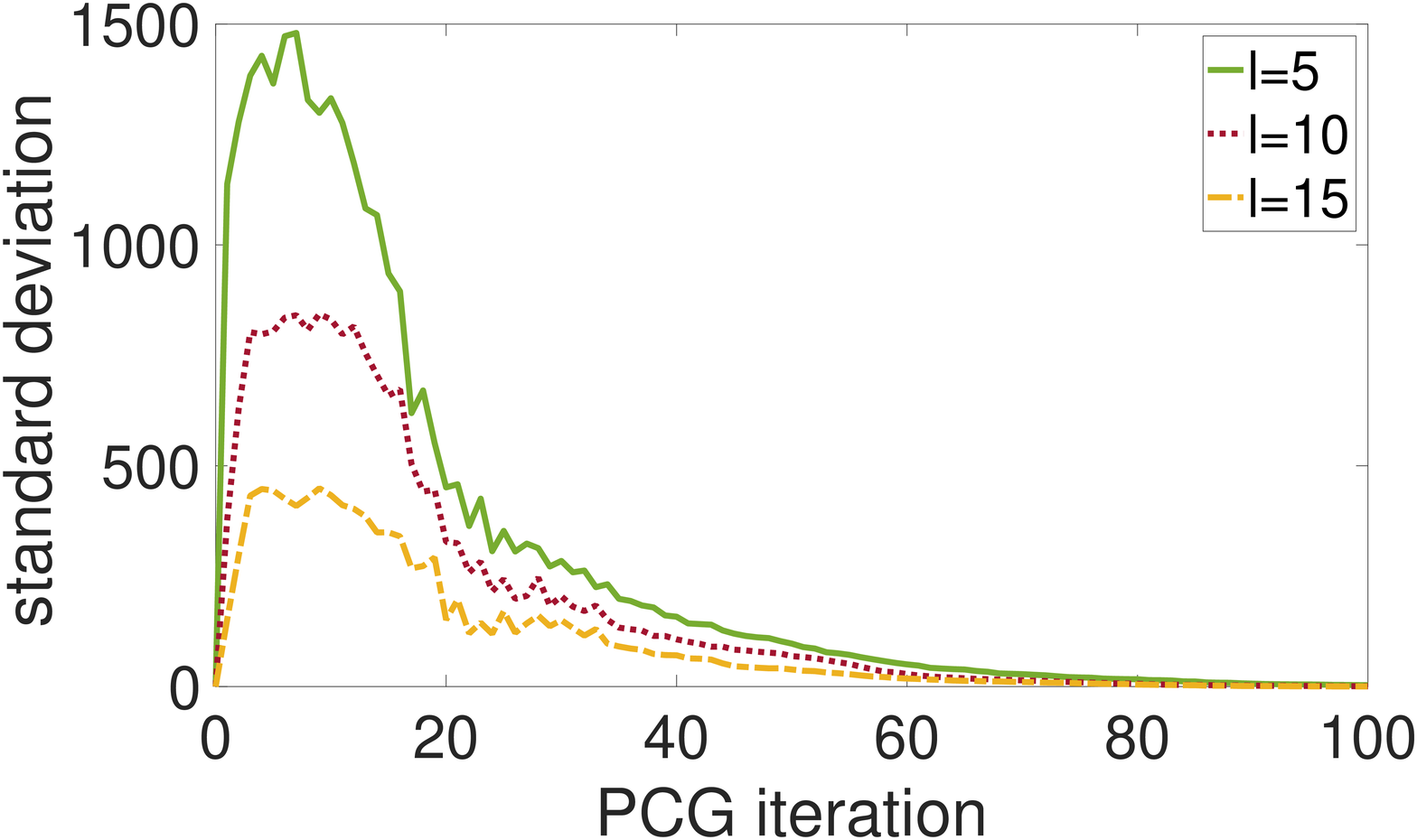}
     \caption{REVD} \label{fig:std_revd_ft10_fx10_ps37_k15_l-5-10-15}
\end{subfigure}
\begin{subfigure}[b]{0.5\linewidth}
  \centering
 \includegraphics[width=\linewidth]{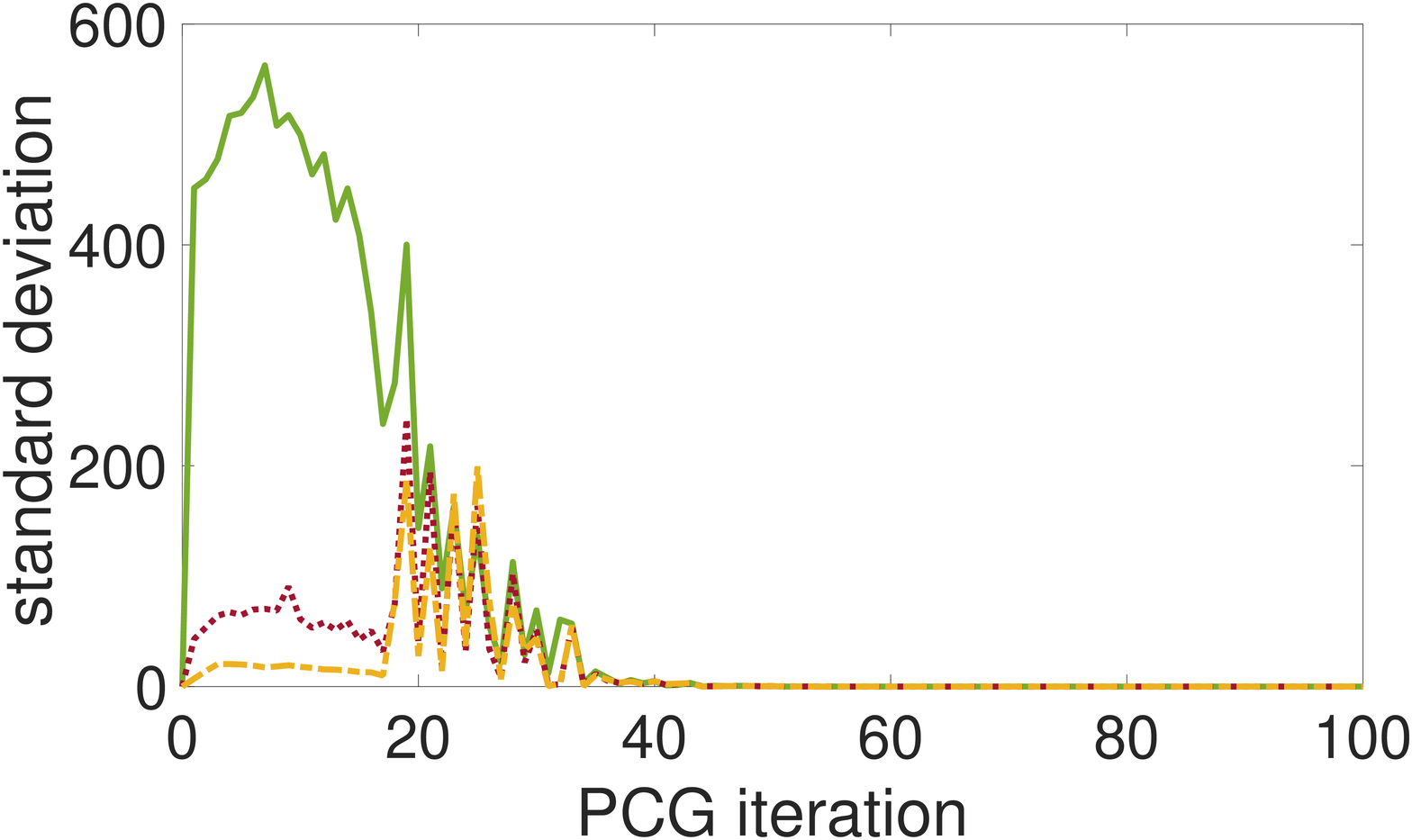}
   \caption{Nystr\"{o}m} \label{fig:std_nystrom_ft10_fx10_ps37_k15_l-5-10-15}
\end{subfigure}\\[1ex]
\begin{subfigure}[b]{0.5\linewidth}
  \centering
 \includegraphics[width=\linewidth]{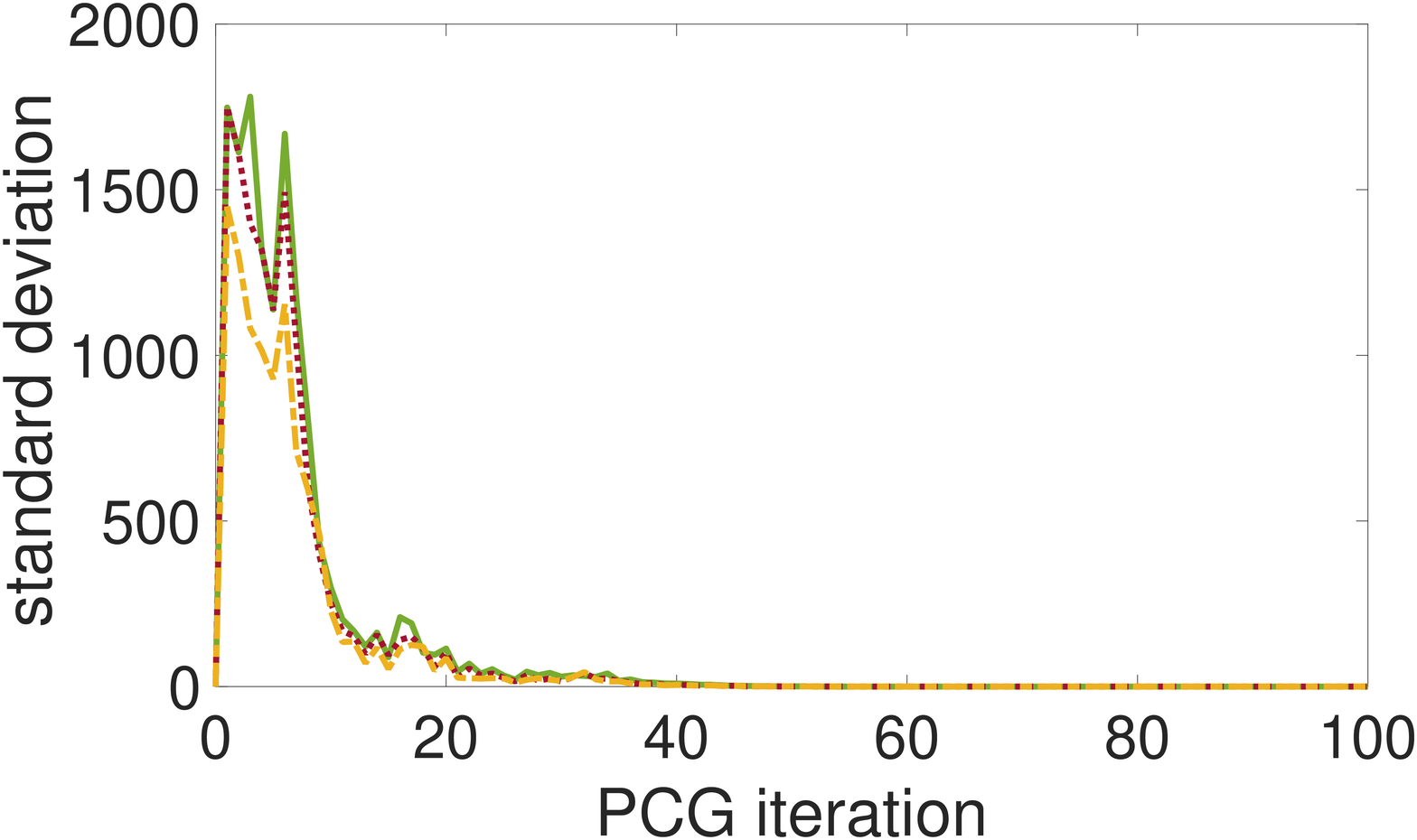}
     \caption{REVD\_\textit{ritzit}} \label{fig:std_ritzit_ft10_fx10_ps37_k15_l-5-10-15}
\end{subfigure}
\caption{Standard deviation of the quadratic cost function at every iteration of PCG when the spectral-LMP is constructed with different randomised methods. For every randomised method we do 50 experiments. Here $\sigma_q = 0.1$ and $L_q = 2 \Delta X$.} 
\label{fig:std_ft10_fx10_ps37_k15_l-5-10-15}
\end{figure}

We next consider the effect of increasing the value of the oversampling parameter $l$. The observation network is as in Section~\ref{sec:numerics_lorenz_cf} (120 observations in total). We set $k=15$ and perform the second inner loop 50 times for every value of $l \in \{5,10,15\}$ with all three randomised methods. The standard deviation of the value of the quadratic cost function at every iteration is presented in Figure~\ref{fig:std_ft10_fx10_ps37_k15_l-5-10-15}.

For all the methods, the standard deviation is greatest in the first iterations of PCG. It is reduced when the value of $l$ is increased and the largest reduction happens in the first iterations. However, REVD\_\textit{ritzit} is the least sensitive to the increase of the oversampling. With all values of $l$, REVD\_\textit{ritzit} has the largest standard deviation in the first few iterations, but it stills gives the largest reduction of the quadratic cost function. Hence, large oversampling is not necessary if REVD\_\textit{ritzit} is used.

\section{Conclusions and future work}
We have proposed a new randomised approach to second level preconditioning of the incremental weak constraint 4D-Var forcing formulation. It can be preconditioned with an LMP that is constructed using approximations of eigenpairs of the Hessian. Previously, by using the Lanzcos and CG connection these approximations were obtained at a very low cost in one inner loop and then used to construct the LMP in the following inner loop. We have considered three methods (REVD, Nystr\"{o}m and REVD\_\textit{ritzit}) that employ randomisation to compute the approximations. These methods can be used to cheaply construct the preconditioner in the current inner loop, with no dependence on the previous inner loop, and are parallelisable.

Numerical experiments with the linear advection and Lorenz-96 models have shown that the randomised LMPs constructed with approximate eigenpairs improve the convergence of PCG more than deterministic LMPs with information from the previous loop. The quadratic cost function reduces more rapidly when using a randomised LMP rather than a deterministic LMP, even if the randomised LMP is constructed with fewer vectors than the deterministic LMP. Also, for the randomised LMPs, the more information about the system we use (i.e. more approximations of eigenpairs are used to construct the preconditioner), the greater the reduction in the quadratic cost function. Using more information to construct a deterministic LMP may not result in larger reduction of the quadratic cost function, especially in the first iterations of PCG, which is in line with results in \cite{Tshimanga08}. However, if not enough information is included in the randomised LMP, then preconditioning may have no effect on the first few iterations of PCG. 

Of the randomised methods considered, the best overall performance was for REVD\_\textit{ritzit}. However, if we run a small number of PCG iterations, the preconditioners obtained with different randomised methods give similar results. The performance was independent of the choice of the random Gaussian start matrix and it may be improved with oversampling. 

In this work we apply randomised methods to generate a preconditioner, which is then used to accelerate the solution of the exact inner loop problem \eqref{eq:Ax_eq_b} with PCG method (as discussed in Section~\ref{sec:randomised_evd}). A different approach has been explored by \cite{Bousserez2018} and \cite{Bousserez2020}, who presented and tested a randomised solution algorithm called the Randomized Incremental Optimal Technique (RIOT) in data assimilation. RIOT is designed to be used instead of PCG and employs a randomised eigenvalue decomposition of the Hessian (using a different method than the ones presented in this paper) to directly construct the solution $\bold{x}$ in \eqref{eq:Ax_eq_b}, which approximates the solution given by PCG.

The randomised preconditioning approach can also be employed to minimise other quadratic cost functions, including the strong constraint 4D-Var formulation. Further exploration of other single-pass versions of the randomised methods for the eigenvalue decomposition, that are discussed in \cite{Halko11}, may be useful. In particular, the single-pass version of the Nystr\"{o}m method is potentially attractive. If a large number of Ritz vectors are used to construct the preconditioner, more attention can be paid to choosing the value of the oversampling parameter $l$ in the randomised methods. In some cases a better approximation may be obtained if $l$ linearly depends on the target rank of the approximation \cite{Nakatsukasa2020}.

\section*{Acknowledgements}
We are grateful to Dr. Adam El-Said for his code for the weak constraint 4D-Var assimilation system. We would like to thank two anonymous reviewers, whose comments helped us to improve the manuscript.

\section*{Conflict of interest}
The authors declare no conflict of interest.

\section*{Funding information}
UK Engineering and Physical Sciences Research Council, Grant/Award Number: EP/L016613/1; European Research Council CUNDA project, Grant/Award Number: 694509; NERC National Centre for Earth Observation.

\bibliography{../../../../PhD_biblio}

\end{document}